\documentclass[onecolumn,final,a4paper]{elsarticle}

\usepackage[utf8]{inputenc}
\usepackage{amsmath,amssymb,amsfonts,latexsym,stmaryrd}
\usepackage{algorithmic}
\usepackage{booktabs}
\usepackage{multirow}
\usepackage[colorlinks=true,breaklinks=true,linkcolor=lightblue,citecolor=lightblue]{hyperref}
\usepackage{stackengine}
\usepackage{mathrsfs}
\usepackage{rotating} 
\usepackage{showkeys}
\usepackage{amssymb,amsmath,bbm,amsfonts,latexsym,subfigure}
\usepackage{amsthm}
\usepackage{graphicx,float,epsfig,color,fancyhdr}
\usepackage{amsopn}

\usepackage[utf8]{inputenc}
\usepackage{graphicx}
\usepackage{wrapfig}

\usepackage{amsfonts}
\usepackage{amsmath}
\numberwithin{equation}{section}
\usepackage{mathtools}
\usepackage{todonotes}
\setuptodonotes{fancyline, color=green!30}

\usepackage[ruled,vlined,linesnumbered]{algorithm2e}
\let\oldnl\nl
\newcommand{\nonl}{\renewcommand{\nl}{\let\nl\oldnl}}%

\usepackage[utf8]{inputenc}
\usepackage{graphicx}
\usepackage{wrapfig}

\usepackage{amsfonts}
\usepackage{amsmath}
\numberwithin{equation}{section}
\usepackage{mathtools}

\usepackage{makecell}

\usepackage[ruled,vlined,linesnumbered]{algorithm2e}
\let\oldnl\nl

\newcommand{\tenprod}{\circ}
\newcommand{\ten}[1]{\mathcal{#1}}
\newcommand{\mat}[1]{\mathbf{#1}}
\renewcommand{\vec}[1]{\mathbf{#1}}

\newcommand{\PGRAPH}[1]{\noindent\textbf{#1}}
\newcommand{\TTf}{TT}
\newcommand{\QTT}[1]{QTT {#1}}
\newcommand{\HAT}[1]{\widetilde{}}

\usepackage{lipsum}

\newtheorem{theorem}{Theorem}[section]

\def\bx{\mathbf{x}}

\newcommand{\vertiii}[1]{{\left\vert\kern-0.25ex\left\vert\kern-0.25ex\left\vert #1 
    \right\vert\kern-0.25ex\right\vert\kern-0.25ex\right\vert}}

\begin{document}

\begin{frontmatter}

\title{Space-Time Spectral Element Tensor Network Approach for Time Dependent Convection Diffusion Reaction Equation with Variable Coefficients}
\author[TDIV]{Dibyendu Adak}
\author[TDIV]{Duc P. Truong}
\author[SNL]{Radoslav Vuchkov}
\author[SNL]{Saibal De}
\author[CCSDIV]{Derek DeSantis}
\author[SNL]{Nathan V. Roberts}
\author[TDIV]{Kim \O. Rasmussen}
\author[TDIV]{Boian S. Alexandrov}

\address[TDIV]{Theoretical Division,
Los Alamos National Laboratory, Los Alamos, NM 87545, USA}
\address[SNL]{Sandia National Laboratories, Albuquerque, NM, USA}
\address[CCSDIV]{Computer, Computational and Statistical Sciences Division,
Los Alamos National Laboratory, Los Alamos, NM 87545, USA}

\begin{abstract}
In this paper, we present a new space-time Petrov-Galerkin-like method. This method utilizes a mixed formulation of Tensor Train (TT) and Quantized Tensor Train (QTT), designed for the spectral element discretization (Q1-SEM) of the time-dependent convection-diffusion-reaction (CDR) equation. We reformulate the assembly process of the spectral element discretized CDR to enhance its compatibility with tensor operations and introduce a low-rank tensor structure for the spectral element operators. Recognizing the banded structure inherent in the spectral element framework's discrete operators, we further exploit the QTT format of the CDR to achieve greater speed and compression. Additionally, we present a comprehensive approach for integrating variable coefficients of CDR into the global discrete operators within the TT/QTT framework. The effectiveness of the proposed method, in terms of memory efficiency and computational complexity, is demonstrated through a series of numerical experiments, including a semi-linear example.
\end{abstract}

\end{frontmatter}


\section{Introduction}
\subsection{Model Time-dependent 3D Convection-Reaction-Diffusion  Problem}
This paper develops a space-time spectral elements method for tensor network low-rank numerical solution of the time-dependent 3D convection diffusion reaction (CDR) equation with Dirichlet boundary conditions of the type,
\begin{equation}
  \begin{split}
  \frac{\partial u}{\partial t}  -\nabla \cdot \Big( {\kappa}(t,\mathbf{x}) \nabla u \Big)+\vec{b} (t,\mathbf{x}) \cdot \nabla u 
   +c(t,\bx)u &= f(t,\bx) \ \text{in} \ [0,T] \times \Omega  ,\\
    u&=g(t,\mathbf{x}) \quad \text{on} \ [0,T] \times  \partial \Omega ,\\
   u(t_0=0,\mathbf{x})&=u_0(\mathbf{x}) \quad \text{in} \ \Omega,
  \end{split}
    \label{model:prob1}
\end{equation}
 where ${\kappa}(t,\mathbf{x})$, $\vec{b} (t,\bold{x}):=[b_1(t,\bold{x}),b_2(t,\bold{x}),b_3(t,\bold{x})]^T$, and $c(t,\bx)$ are the nonlinear coefficients, and $\bx=(x,y,z)$. The computational domain in space $\Omega = \Omega_Z \times \Omega_Y \times \Omega_X  \subset \mathbb{R}^3$  is a three dimensional cube which is a Cartesian product of three intervals, and $T$ is the final time-point. We have considered inhomogeneous boundary condition and initially the tensor-train space-time formats are developed for coefficients that are constants or allow separation of variables, and later extended for a general type of variable coefficients. 
\subsection{Classical Methods for Solving CDR}
Spectral element method (SEM) is a powerful technique for solving CDR with high accuracy \cite{sengupta2020global}. By combining the geometric flexibility of finite elements with the exponential convergence of spectral methods, SEM can achieve excellent resolution of complex solutions \cite{codina1998comparison}. 
Discontinuous  Galerkin \cite{shu2009discontinuous} and Petrov-Galerkin \cite{bui2013unified} methods are  widely used for solving CDR, particularly in finite element and spectral methods \cite{gharibi2021convergence}. These methods involve approximating the PDE solution  by projecting it onto a finite-dimensional subspace of trial functions. The key idea is to ensure that the residual error is orthogonal to the chosen subspace of test functions. 

In the Galerkin approach, the continuous problem is reformulated into a discrete problem by selecting a finite number of trial functions to approximate the solution. These functions are typically chosen from a function space like polynomials or piecewise-defined functions. In the classical (Bubnov-)Galerkin method, the same set of functions is used for both the trial (approximation of the solution) and the test (weighting) functions, which ensures that the residual error is orthogonal to the space spanned by the trial functions. The PDE is transformed into its weak (integral) form by multiplying it by a test function and integrating over the domain. The weak form allows handling problems with lower regularity solutions and dealing with complex geometries. 

The Petrov-Galerkin (PG) method generalizes the Bubnov-Galerkin method by allowing test and trial spaces to differ. This added flexibility allows for improved stability and accuracy, especially for convection-dominated problems or where numerical instability is a concern. By selecting appropriate test functions, the Petrov-Galerkin method can introduce stabilization mechanisms, such as upwinding in the case of convection-dominated problems and reducing oscillations or instabilities. In the present work, we do not apply such stabilization mechanisms since our underlying problem is a diffusion-dominated problem. In fact, we employ the \emph{same} discrete basis for both test and trial functions.  Our method is following the Petrov-Galerkin framework because the underlying continuous function spaces differ, specifically in the norms employed; our stability analysis depends on this difference, and it is based on \emph{inf-sup} condition \cite{boffi2013mixed,steinbach2007numerical}.

Space-time discretizations treat the spatial and temporal domains simultaneously with the same methodology, whereas time-marching schemes discretize the spatial domain with one methodology e.g., FEM, Nonconforming FEM, and recently introduced virtual element method(VEM), and employ a distinct time-stepping scheme (such as backward Euler) for the temporal discretization.  Space-time methods usually demand more resources, especially in the memory required for the representation of the solution across a time slab, but they allow more flexibility in the temporal discretization, including the possibility of spatially-localized refinement in time.  Space-time methods can also be easier to analyze and may have higher accuracy than a corresponding time-marching scheme. Our current numerical examples serve as a proof of concept for the new tensor method. In future work, we aim to tackle more numerically challenging scenarios where stability may be a concern.


\subsection{Mitigating the Curse of Dimensionality}
 The CDR equations model a wide range of phenomena in physics, chemistry, and engineering, which often demand substantial computational resources due to the need for fine spatial and temporal discretizations. Traditional approaches to solving such problems often lead to prohibitively large systems of equations, making them impractical for many real-world applications, since the number of grid points grows exponentially with the number of dimensions. In real-life applications modeled by time-dependent PDEs, such as full waveform inversion problems, the number of grid points required for solving the PDE can be as large as $M = 6.6 \times 10^{10}$ per time step. With a large number of time steps  $N = 4 \times 10^5$ , one must store a total of  $MN$  floating point numbers \cite{krebs2009fast}.
 
 This phenomenon is known as the \emph{curse of dimensionality} \cite{bellman1966dynamic}. The curse of dimensionality leads to poor computational scaling in numerical algorithms and poses the main challenge in multidimensional numerical computations, irrespective of the specific problem. Notably, even exascale high-performance computing, with its optimization strategies, cannot overcome the curse of dimensionality. This phenomenon, forces algorithm developers to make hard choices—either to reduce the fidelity of the model using, e.g., reduced-order models (ROMs) or to repeat computation using checkpointing type strategies. ROMs focus on creating a reduced form of the PDE operator and are an excellent solution in some cases. However, they are often tailored for specific models and demand significant domain expertise. Furthermore, ROMs can be incompatible with legacy codes as they often require invasive modification of the simulation software. ROM development is an active area of research and can work well for solving certain classes of problems \cite{leugering2014trends,gunzburger2007reduced}.

As an alternative, tensor network techniques have shown promise in alleviating the curse of dimensionality associated with high-dimensional problems. Among these, the tensor train (TT) format, introduced by Oseledets and Tyrtyshnikov \cite{oseledets2010tt}, has gained particular attention due to its ability to represent high-dimensional data efficiently while maintaining computational tractability.

Here we  present a novel approach that combines the high accuracy of space-time spectral element methods with the computational efficiency of TT decomposition for solving time-dependent CDR equations in three spatial dimensions. Our method exploits the inherent low-rank structure often present in the solutions of such equations to dramatically reduce the computational complexity and storage requirements.
Moreover, given that the spectral element discrete operators have banded structures, we further exploit the Quantized Tensor Train (QTT) format for more economical and efficient solvers~\cite{kazeev2013multilevel}.
Specifically, we develop a mixed TT/QTT decomposition of the four-dimensional (three space plus time) space-time spectral element Petrov-Galerkin discretization, and demonstrate how to construct the TT and QTT representations of the spectral element discretization.
Our results show that this TT/QTT-SEM-PG approach can achieve high accuracy while dramatically reducing the computational resources required, making it possible to solve previously intractable problems in CDR systems.

The remainder of this paper is organized as follows. In Section~\ref{sec:SEM model}, we derive the weak formulation of \eqref{model:prob1}, and  review some basic concepts for space-time spectral element methods, and introduce the discretization and its matrix formulation. Section~\ref{sec:aseemble} describes the assembly process that leads to the global linear system. We have reformulated this process to be more compatible with tensor operations. This reformulation is particularly important as it simplifies and optimizes the tensorization process.
%
%
In Section~\ref{sec:Tensor-Networks}, we present our mixed TT/QTT design of the
numerical solution of the CDR equation and introduce our algorithms in TT/QTT format. In
Section~\ref{sec:numerical}, we present our numerical results and assess the performance of our method. In Section~\ref{sec:conclusion}, we offer our final remarks and discusses possible future work. Additional details on the CDR tensorization approach are provided in the appendices.
\section{Space-Time Finite Element Formulation of CDR}
\label{sec:SEM model}
The Spectral Element Method (SEM) is a high-order numerical method used to solve partial differential equations (PDEs). Seen as a hybrid between the spectral method and the finite element method (FEM), it combines the geometric flexibility of FEM with the high accuracy of spectral methods, making it particularly well-suited for problems involving complex geometries or requiring high precision. In this section, we provide the background details of the SEM used throughout the text. 

We will employ the following notations:  $\Omega_X,\Omega_Y,\Omega_Z \subset \mathbb{R}$ will denote fixed intervals of $\mathbb{R}$, while   $\Omega = \Omega_X \times \Omega_Y \times \Omega_Z \subset \mathbb{R}^3$. The time interval will be denoted by $I_T=[0,T]$, and our space-time domain $\Omega_T = [0,T] \times \Omega$.


In FEM, the domain is broken into smaller, simple ``elements'' and we approximate the solution $u$ using simple functions on the elements. This is done by first transforming the PDE into its weak formulation.  The domain is divided into elements for defining basis functions.  The weak formulation is then discretized in terms of the basis functions, from which the mass and stiffness matrices are computed. The resulting linear system is then formulated and solved.  Space-time FEM is formulated over the entire space-time domain $\Omega_T$.  In this subsection, we detail the FEM for the 3D CRD in Equation (\ref{model:prob1}).

\subsection{Space-Time Weak Formulation}

We let $C_c^\infty(\Omega)$ denote the set of infinitely differentiable functions with compact support within open set $\Omega$, while $L^2(\Omega)$ and $L^\infty(\Omega)$ will denote the set of (equivalence classes of) square integrable and essentially bounded functions over $\Omega$ respectively.    Given a multi-index $\alpha = (\alpha_1, \alpha_2,\alpha_3)$ consisting of non-negative integers, and a function $\phi \in C_c^\infty(\Omega)$, we let $D^{\alpha}\phi$ denote the $\alpha$\textsuperscript{th} partial derivative of $\phi$:
\[
D^{\alpha}\phi = \frac{\partial^{|\alpha|} \phi}{\partial x^{\alpha_1}\partial y^{\alpha_2}\partial z^{\alpha_3}}
\]
where $|\alpha| = \alpha_1 + \alpha_2 + \alpha_3$. If $u$ is a (locally) integrable function, then we let $D^{\alpha}u$ denote the $\alpha$\textsuperscript{th} partial derivative of $u$ in the weak sense.  That is, if there exists a (locally integrable) function $v$ such that
\[
\int_\Omega u D^{\alpha}\phi = (-1)^{|\alpha|} \int_\Omega v \phi
\]
for all $\phi \in C_c^\infty(\Omega)$, then we assign $D^\alpha u = v$. We let $H^1(\Omega)$ denote the Sobolev space consisting of all square-integrable functions which are differentiable in the weak sense up to first order. We equip $H^1(\Omega)$ with with norm
\[
\lVert u \rVert_{H^1(\Omega)}:= \Big ( \sum_{|\alpha| \leq 1} \lVert D^{\alpha} u \lVert^2  \Big)^{1/2},
\]
where $\lVert\cdot\rVert$ is the standard $L^2(\Omega)$ norm. The space $H_0^1(\Omega)$ denotes the special subspace of $H^1(\Omega)$ consisting of functions that vanish on the boundary in the weak sense, and thus satisfy homogeneous boundary conditions.  
By $H^{-1}(\Omega)$, we mean the dual space of $H_0^1(\Omega)$, i.e. the space of all bounded linear functionals on  $H_0^1(\Omega)$.  
We give $H^{-1}(\Omega)$ the standard operator norm:
\[
\lVert u\rVert_{H^{-1}(\Omega)} = \operatorname*{sup}_{\phi \in H^{1}_0(\Omega), \|\phi\|_{H^1(\Omega)} \leq 1} | \langle u, \phi \rangle |.
\]
Given a function space $X$ with associated norm $\lVert\cdot\rVert_X$, the Bochner space $L^2(0, T; X)$ consists of (equivalence classes of) functions $u$ such that $u(t, \cdot) \in X$ for almost all $t \in [0, T]$. We equip this space with the norm
\begin{equation*}
    \lVert u \rVert_{L^2(0, T; X)} = \left(\int_0^T \lVert u(t, \cdot) \rVert_X^2 dt \right)^{\frac{1}{2}}.
\end{equation*}
Similarly, we define the space $H^1(0, T; X)$ as the space of all functions whose weak time derivatives are in $L^2(0, T; X)$; we equip it with the norm
\begin{equation*}
    \lVert u \rVert_{H^1(0, T; X)} = \left(\lVert u \rVert^2_{L^2(0, T; X)} + \left\lVert \frac{\partial u}{\partial t} \right\rVert^2_{L^2(0, T; X)}\right)^{\frac{1}{2}}.
\end{equation*}
These spaces, especially with $X = H^1(\Omega)$ and $X = H^{-1}(\Omega)$, define the appropriate function spaces for posing the weak formulation of CRD.
The space $L^2(0, T; H^1(\Omega))$ ensures a function $u(t, \bold{x})$ is square-integrable in time, and each time slice $u(\cdot, \bold{x})$ has square-integrable first spatial derivatives. Similarly, the space $H^1(0, T; H^{-1}(\Omega))$ ensures that a function $u(t, \bold{x})$ has first order time derivative that is weakly square integrable in time.

We define the trial space $U := L^2(0,T;H^{1}_0(\Omega)) \cap H^1(0, T; H^{-1}(\Omega))$ associated with the norm $\|u\|^{2}_U = \|u\|^{2}_{L^2(0,T;H^{1}_0(\Omega))} + \|u\|^{2}_{H^1(0, T; H^{-1}(\Omega))}$, and the test space $V:=L^2(0,T;H_0^{1}(\Omega))$. 
Integrating by parts, we obtain the weak formulation of the model problem \eqref{model:prob1} with homogeneous boundary condition as finding the function $u \in U$ such that
\begin{equation}
 \left\langle  \frac{\partial u}{\partial t}, v \right\rangle + \langle \kappa(t,\bold{x}) \nabla u , \nabla v \rangle + \langle \vec{b} (t,\bold{x}) \cdot \nabla u, v \rangle  + \langle  c(t,\bx)u, v \rangle  =   \langle  f(t,\bx),v \rangle \quad \forall v \in V.
\label{Weak:Form}
\end{equation}
Here, we assume that the forcing function is in dual space of $H^1_0(\Omega)$ for a.e. $t \in I_T$, namely  $f(t, \textbf{x}) \in L^2(0, T; H^{-1}(\Omega))$, and the coefficients $\kappa(t, \textbf{x}), b_i(t, \textbf{x}), c(t, \textbf{x}) \in L^{\infty}(0,T; L^{\infty}(\Omega))$.
We define the bilinear term
\begin{equation}
\label{Biinear:FM}
    \mathcal{D}(u,v):= \left\langle  \frac{\partial u}{\partial t}, v \right\rangle + \langle \kappa(t,\bold{x}) \nabla u , \nabla v \rangle + \langle \vec{b} (t,\bold{x}) \cdot \nabla u, v \rangle  + \langle  c(t,\bx)u, v \rangle.
\end{equation}
Further, by assuming the we state the boundedness of the bilinear form $\mathcal{D}(u,v)$. 
\begin{equation}
\label{bound:D}
    |\mathcal{D}(u,v) | \leq C(\kappa^*, \| \bold{b}\|_{\infty,\Omega_T}, \| c\|_{\infty,\Omega_T}) \| u\|_{U} \| v\|_{V}, 
\end{equation}
where, $\kappa^*=\text{ess sup} \{ \kappa(t,\bold{x}) : (t,\bold{x}) \in [0,T] \times \Omega \}$.
Moreover the bilinear from $\mathcal{D}(\cdot,\cdot)$ satisfies the following stability condition \cite[Theorem~2.1]{steinbach2015space} 
\begin{equation}
\label{stability:D}
    C \|u\|_{U} \leq  \sup_{0 \neq v \in V} \frac{\mathcal{D}(u,v)}{\|v \|_{V}} \quad \forall v \in V.
\end{equation}
We are in a position to state unique solvability of the Petrov-Galerkin variational formulation~\eqref{Weak:Form}:
\begin{theorem}
\label{wellposed:weak}
    Let us assume that $f \in L^2(0,T; H^{-1}(\Omega))$, and  the bilinear form $\mathcal{D}(\cdot,\cdot)$ \eqref{Biinear:FM} is bounded as discussed in \eqref{bound:D} and satisfies stability condition \eqref{stability:D}. Then, there exists a unique solution $u \in U$ of \eqref{Weak:Form} satisfying 
    \begin{equation}
        \|u\|_{U} \leq C \| f\|_{L^2(0,T; H^{-1}(\Omega))},
    \end{equation}
    where $C$ is positive constant independent of mesh size $h$ but depends on $\kappa(t,\bold{x})$, $\bold{b}(t,\bold{x})$, and $c(t,\bold{x})$.
\end{theorem} 
Interested readers can refer \cite[Theorem~3.7]{boffi2013mixed,steinbach2007numerical} for a detailed proof of Theorem~\ref{wellposed:weak}. 
In the weak formulation of the CDR in Equation \eqref{Weak:Form}, the solution $u$ and test functions $v$ are chosen from appropriate function spaces to account for both the spatial and temporal regularity. Contrast this against the (Bubnov-)Galerkin methods, where we pick the same function space for the trial and test functions, which may cause instability. For example, when the velocity field is strong (i.e., large convection), the convective term can cause steep gradients or sharp layers to appear in the solution, which the standard Galerkin method fails to resolve properly. In the our Petrov-Galerkin approach, the test functions $v$ can come from a different space than the trial functions.
Equation \eqref{Weak:Form} is  the weak formulation corresponding to \eqref{model:prob1} with homogeneous boundary conditions \cite{larson2010finite,steinbach2015space,gomez2024space}. The weak formulation corresponding to inhomogeneous boundary condition can be obtained following the standard technique of decomposing the solution, $u(t,\bold{x})$, of \eqref{model:prob1} to: $u(t,\bold{x})= u^{\text{homog}}(t,\bold{x})+u_0(t,\bold{x})$, where $u^{\text{homog}}(t,\bold{x})$ is the unknown solution of \eqref{model:prob1} with homogeneous Dirichlet boundary conditions, and $u_0(t,\bold{x})$ satisfies the boundary conditions at $[0,T] \times \Omega$, and it is different from zero only on the boundaries. 
In the weak formulation this simply results in a modified loading term, $f(t, \bold{x})$.

%

\subsection{Finite Element Approximation}

We now describe the finite element approximation of the Equation \eqref{Weak:Form}.  This is done by picking finite dimensional subspaces $U_h \subset U$ and $V_h \subset V$ specified by chosen bases described below.  Let $\{\phi_i\}_{i=1}^N$ be a basis for $U_h$ and $\{\psi_j\}_{j=1}^M$ be a basis for $V_h$. Then any approximate solution $u_h \in U_h$ can be written as 
\[
u_h(t, \textbf{x}) = \sum_{i=1}^N U_i \phi_i(t, \textbf{x})
\]
and each test function $v_h \in V_h$ is written as
\[
v_h(t, \textbf{x}) = \sum_{j=1}^M V_j \psi_j(t, \textbf{x}).
\]
Substituting these into the weak form in Equation \eqref{Weak:Form} results in a discrete system.  We now describe this system for a particular choice of linear basis functions.

First, the domain $\Omega_T$ is discretized into a set of hypercubes which overlap only on their boundaries.  The hypercube elements are defined by dividing each axis into smaller intervals: 
\begin{align}
    \Omega_X &= [x_0, x_1] \cup [x_1,x_2] \cup \dots [x_{n_X-1}, x_{n_X}]\notag \\
    \Omega_Y &= [y_0, y_1] \cup [y_1,y_2] \cup \dots [y_{n_Y-1}, y_{n_Y}]\notag  \\
    \Omega_Z &= [z_0, z_1] \cup [z_1,z_2] \cup \dots [z_{n_Z-1}, x_{n_Z}]\notag  \\
    [0,T] &= [t_0, t_1] \cup [t_1,t_2] \cup \dots [t_{n_T-1}, t_{n_T}] \notag 
\end{align}
The resulting hypercubes are of the form $[x_i,x_{i+1}] \times [y_j,y_{j+1}] \times [z_k, z_{k+1}] \times [t_l, t_{l+1}]$. Throughout, we adopt a uniform mesh on $\Omega_T$ with $h:=x_{i+1}-x_i =y_{i+1}-y_i =z_{i+1}-z_i =t_{i+1}-t_i$ for all $i,j,k,l$.  Further, we assume $N=n_X=n_Y=n_z=n_T$ so that there is a total of $N+1$ total nodes along each dimension, and $N$ intervals. We let $k=0,1,\dots N$ denote the index for the nodes in each of the one-dimensional spaces.  There is a total of $N_Q = (N+1)^4$ total space-time elements in $\Omega_T$.   We denote these elements by $\{q_l\}_{l=0}^{N_Q-1}$, which are raster ordered with $X$ before $Y$ before $Z$ before $T$. The bijection between the linear element index $l$ and the individual element indices $(k_x, k_y, k_z, k_t)$ is then given by
\[
l \longleftrightarrow k_x+(N+1)k_y+(N+1)^2k_z+(K+1)^3k_t.
\]

\subsubsection{Local Estimates}

Throughout we make use of local interpolation, which we describe here. Each 4D element $q_l$ consists of two nodes per dimension for a total of 16 total nodes. For a fixed element, each node can be indexed by four indices $i_1, i_2, i_3, i_4 \in \{0,1\}$, with the zero index corresponding to the lower value in the interval.  The \emph{local index} for each element $i=0,1,\dots, 15$ is given by
\[
i = i_1 + 2 i_2 + 4i_3 + 8 i_4
\]
Along each dimension, we define the standard piecewise linear hat functions defined on the one-dimensional elements.  For example on the interval $I_k^x := [x_{k},x_{k+1}]$
\begin{equation}
\label{eq: phi}
\phi^x_{0,k}(x):= \begin{cases} 
      \frac{x - x_{k}}{h} & x \in I_k^x  \\
      0 & \mbox{else}
   \end{cases} \quad ~\text{and}~ \quad 
\phi^x_{1,k}(x) := \begin{cases} 
      \frac{x_{k+1}-x}{h} & x \in I_k^x \\
      0 & \mbox{else}.
   \end{cases} 
\end{equation}
Clearly, the $\phi^x_{i_1,k}$ have almost everywhere defined derivative
\begin{equation}
\label{eq: derivative phi}
\frac{\partial \phi^x_{0,k}}{\partial x}(x):= \begin{cases} 
      \frac{1}{h} & x \in (x_{k-1},x_k)  \\
      0 & \mbox{else}
   \end{cases} \quad ~\text{and}~ \quad 
\frac{\partial \phi^x_{1,k}}{\partial x}(x) := \begin{cases} 
      \frac{-1}{h} & x \in (x_{k-1},x_k) \\
      0 & \mbox{else}.
   \end{cases} 
\end{equation}
We similarly define $\phi^y_{0,k},\phi^y_{1,k}, \phi^z_{0,k},\phi^z_{1,k}, \phi^t_{0,k},\phi^t_{1,k}$.
Using the local index notation associated with element $q_l$, the local space-time basis functions defined as the tensor product of the 1D hat functions:
\[
\phi_{i,l}(t,x,y,z) = \phi_{i_1,k_x}(x) \phi_{i_2,k_y}(y) \phi_{i_3,k_z}(z) \phi_{i_4,k_t}(t)
\]
where we have used the correspondence between the element $l$ and interval index $(k_x,k_y,k_z,k_t)$, and the local index $i$ with node index $(i_1,\dots, i_4)$.  For completeness, we note that the 16 functions $\phi_{i,l}$ are only non-zero on the cell $q_l$.


\subsection{Local and Global Interpolation Operators}

Now we will define the global and local Lagrange interpolation operators $L_h$ and $L_h^l$ respectively. For cell $q_l$ and local index $i$, we let $V_{l,i}$ denote the $(t,x,y,z)$ node at index $i$ in cell $l$. We define the local interpolation operator at $q_l$, $L_h^l:U \rightarrow U_h$, as 
\begin{equation}
   u^l_h(\textbf{x},t):=L_h^l(u)(\textbf{x},t)=\sum_{i=0}^{15} u(V_{l,i}) \phi_{i,l}(\textbf{x},t)
    \label{Lagrange:operator}
\end{equation}
 In other words, $L_h^l(u)(\textbf{x},t)$ approximates the solution locally using the hat functions defined on the vertices of the finite element with index $l$. The global interpolation operator $L_h:U \rightarrow U_h$ is defined as
\begin{equation}
    u_h(t,\textbf{x}) := L_h(u)(t,\textbf{x}) = \sum_{l=0}^{N_Q-1} L_h^l(u)(\textbf{x},t) =  \sum_{l=0}^{N_Q-1} \sum_{i=0}^{15} u(V_{l,i}) \phi_{i,l}(\textbf{x},t)
    \label{Dis_uh}
\end{equation}
The global interpolation operator combines the contributions from all elements of the mesh.  These operators may be applied to the solution $u$, as well as the coefficient and forcing functions. 

Let $V_h$ denote the finite dimensional subspace spanned by a set of linearly independent functions $\{\psi_j\}_{j=1}^M \subset V$. Substituting the approximate solution $u_h$ and a $v_h \in V_h$ into the weak form Equation \ref{Weak:Form} leads to our discrete system of equations.  Specifically, the Petrov-Galerkin discretization of the variational Problem  is find $u_h \in U_h$ such that
\begin{equation}
\begin{split}
\underbrace{\left \langle  \frac{\partial u_h}{\partial t}, v_h \right \rangle}_{\bold{T_1}} + \underbrace{\langle L_h(\kappa(t,\bold{x})) \nabla u_h , \nabla v_h \rangle}_{\bold{T_2}} + \underbrace{\langle L_h(\vec{b} (t,\bold{x})) \cdot \nabla u_h, v_h \rangle}_{\bold{T_3}} & + \langle  \underbrace{L_h(c(t,\bx))u_h, v_h \rangle}_{\bold{T_4}} \\&  =   \underbrace{\langle  L_h(f(t,\bx)),v_h \rangle}_{\bold{T_5}},
 \quad \forall v_h \in V_h.
 \end{split}
\label{Weak:Form2}
\end{equation} 
While the $\psi_j$ can have looser regularity than the $\phi_i$ basis for $U_h$, in practice we will consecutively use each of the $2^4$ local space-time Lagrange basis functions, $\phi_i(t,\bold{x})$, as the test function, $v_h$, in Eq.~\eqref{Weak:Form2}.
We denoted each term of Eq. (\ref{Weak:Form2}) with $\bold{T_i}$. 
The goal of this paper is to provide a tensor representation for all   $\bold{T_i}$ of Equation \ref{Weak:Form2} for fast and efficient solving. In the next section, we provide the necessary background on tensor train methods.
Further, we conclude this section with the discrete stability condition of the discrete scheme \eqref{Weak:Form2} \cite{steinbach2015space}.
\begin{theorem}
    Let $U_h \subset U$, and $V_h \subset V$ be the discrete space satisfying $U_h \subset V_h$. Then there holds a discrete stability condition 
    \begin{equation*}
        C \|u_h\|_{U} \leq \sup_{0 \neq v_h \in V_h} \frac{\mathcal{D}(u_h,v_h)}{\|v_h \|_{V}} \quad \forall v_h \in V_h,
    \end{equation*}
    where $C$ is positive constant depend on the regularity of $\kappa(t,\bold{x})$, $\bold{b}(t,\bold{x})$, and $c(t,\bold{x})$, but independent of mesh size $h$.
\end{theorem}
Since the discrete bilinear form is bounded and satisfies the discrete stability condition, the discrete scheme~\eqref{Weak:Form2} is well-posed, i.e., \eqref{Weak:Form2}
has unique solution (see~ \cite[section~3]{steinbach2015space}).



\section{Tensor Train Decomposition}
\label{sec:TT}

In this section, we introduce the \TTf{} format \cite{oseledets2011tensor}, as well as the representation of linear operators in the so-called \textit{TT-matrix} format, the cross-interpolation method, and the QTT format.
All these methods are fundamental in the tensorization of our spectral element discretization of the CDR equation.

\subsection{Tensor Train}
\label{subsec:TNs}
\begin{figure}
    \centering
    \includegraphics[width=0.6\linewidth]{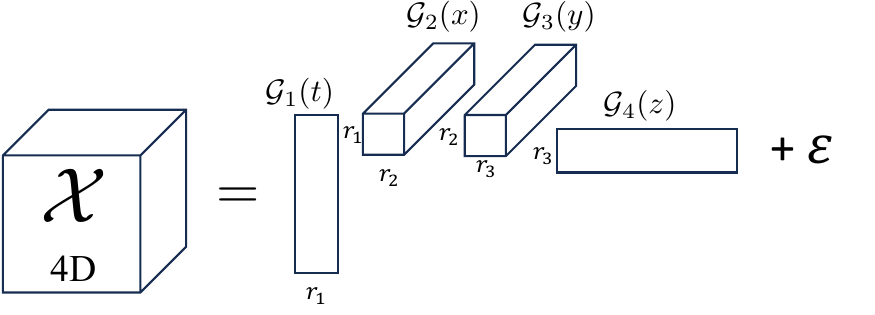}
    \caption{TT format of a 4D tensor with \TTf{} ranks
    $\mathbf{r} = \big[r_1,r_2,r_3\big]$ and approximation error
    $\varepsilon$, in accordance with Eq.~\eqref{eqn:TT_def_element}.}
    \label{fig:TT_4D}
\end{figure}

The TT format, introduced by Oseledets in
2011~\cite{oseledets2011tensor}, represents a sequential chain of
matrix products involving both two-dimensional matrices and
three-dimensional tensors, referred to as TT-cores.
We can visualize this chain as in Figure \ref{fig:TT_4D}.
Given that tensors in our formulation are at most four dimensional
(one temporal and three spatial dimensions), we consider the tensor
train format in the context of 4D tensors.
Specifically, the \TTf{} approximation $\ten{X}^{TT}$ of a
four-dimensional tensor $\ten{X}$ is a tensor with elements
\begin{align}
  \ten{X}^{TT}(i_1,i_2,i_3,i_4)
  =  \sum^{r_1}_{\alpha_{1}=1}\sum^{r_2}_{\alpha_{2}=1}\sum^{r_3}_{\alpha_{3}=1}
  \ten{G}_1(1,i_1,\alpha_1)\ten{G}_2(\alpha_1,i_2,\alpha_2)\ten{G}_3(\alpha_2,i_3,\alpha_3)\ten{G}_4(\alpha_{3},i_4,1).
  \label{eqn:TT_def_element}
\end{align}
Here, we have $\ten{X} = \ten{X}^{TT} + \varepsilon$
where the error, $\varepsilon$, is a tensor with the same dimensions
as $\ten{X}$. The elements of the array $\mat{r} = [r_1,r_2,r_3]$
are the TT-ranks, which quantify the compression effectiveness of the TT approximation.
Since, each \TTf{}-core, $\ten{G}_p(i_k)$, only depends on a single index
of the full tensor $\ten{X}$,  the \TTf{} format effectively
embodies a discrete separation of variables \cite{bachmayr2016tensor}.
In Figure~\ref{fig:TT_4D}, we show a four-dimensional array
$\ten{X}(t,x,y,z)$, decomposed in TT-format.

\subsection{Linear Operators in TT-matrix format}
\label{SUB:Lin_operinTT}
\begin{figure}[h]
  \centering
  \includegraphics[width=0.8\textwidth]{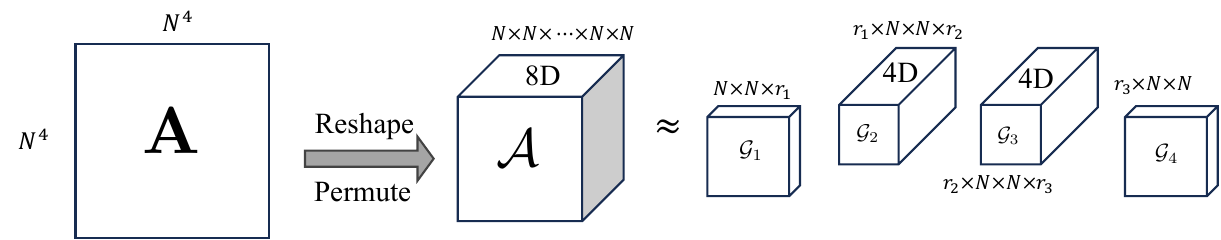}
  \caption{Representation of a linear matrix $\mat{A}$ in the
    TT-matrix format. First, we reshape the operation matrix
    $\bold{A}$ and permute its indices to create the tensor
    $\ten{A}$. Then, we factorize the tensor in the tensor-train
    matrix format according to Eq.~\eqref{eqn:TT-matrix-componentwise}
    to obtain $\ten{A}^{\TTf{}}$.}
  \label{fig:TT-matrix-format}
\end{figure}

Suppose that the approximate solution of the CDR equation is a 4D
tensor $\ten{U}$, then the linear operator $\ten{A}$ acting on that
solution is represented as an 8D tensor.
The transformation $\ten{A}\ten{U}$ is defined as:
\begin{align*}
  \big(\ten{A}\ten{U}\big)(i_1,i_2,i_3,i_4)
  = \sum_{j_1,j_2,j_3,j_4} \ten{A}(i_1,j_1,\dots,i_4,j_4)\ten{U}(j_1,\dots,j_4).
\end{align*}
The tensor $\ten{A}$ can be related to a matrix operator $\bold{A}$ via
\begin{equation}
  \ten{A}(i_1,j_1,\dots,i_4,j_4) = \bold{A}(i_1i_2i_3i_4,j_1j_2j_3j_4).
  \label{eqn:ten_mat_op_relation}
\end{equation}
We can construct the tensor $\ten{A}$ by suitably reshaping and
permuting the dimensions of the matrix $\bold{A}$.
The linear operator $\ten{A}$ can be further represented in a
variant of \TTf{} format, called \emph{TT-matrix},
cf.~\cite{truong2023tensor}.
The component-wise \TTf{}-matrix $\ten{A}^{TT}$ is defined as:
\begin{equation}
  \ten{A}^{TT}(i_1,j_1,\ldots,i_4,j_4)
  =  \sum_{\alpha_{1},\alpha_2,\alpha_{3}}
  \ten{G}_1\big(1,(i_1,j_1),\alpha_1\big)
   \ldots\ten{G}_4\big(\alpha_{3},(i_4,j_4),1\big),
  \label{eqn:TT-matrix-componentwise}
\end{equation}
where $\ten{G}_{k}$ are 4D TT-cores .
Figure~\ref{fig:TT-matrix-format} shows the process of transforming a
matrix operator $\bold{A}$ to its tensor format, $\ten{A}$, and,
finally, to its TT-matrix format, $\ten{A}^{\TTf{}}$.

The Kronecker product $\otimes$ of two matrices is an operation that produces a larger matrix, while the tensor product $\tenprod$ produces a higher dimensional tensor (See Appendix). We can further simplify the TT-matrix representations of the matrix
$\bold{A}$ if it is a Kronecker product of matrices, i.e.
$\bold{A}=\bold{A}_1\otimes\bold{A}_2\otimes\bold{A}_3\otimes\bold{A}_4$.
Based on the relationship defined in
Eq.~\eqref{eqn:ten_mat_op_relation}, the tensor $\ten{A}$ can be
constructed using tensor product as $\ten{A} = \bold{A}_1 \tenprod
\bold{A}_2 \tenprod \bold{A}_3 \tenprod \bold{A}_4$.
This implies the internal ranks of the TT-format of $\ten{A}$
in \eqref{eqn:TT-matrix-componentwise} are all equal to $1$.
In such a case, all summations in
Eq.~\eqref{eqn:TT-matrix-componentwise} reduce to a sequence of single
matrix-matrix multiplications, and the \TTf{} format of $\ten{A}$ becomes
the tensor product of $d$ matrices:
\begin{align}
  \ten{A}^{TT} &= \mat{A}_1 \tenprod \mat{A}_2 \tenprod \dots \tenprod \mat{A}_d.
  \label{eqn:TT_matrices}
\end{align}
This specific structure appears quite often in the matrix
discretization, and will be exploited in the tensorization to
construct efficient \TTf{} format.

\subsection{TT Cross Interpolation}
\label{sec:Cross Interpolation}
\begin{figure}
    \centering
    \includegraphics[width=0.5\linewidth]{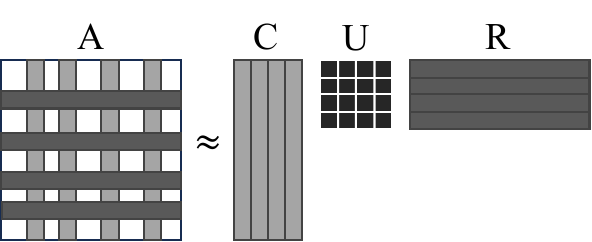}
    \caption{CUR matrix decomposition.}
    \label{fig:CUR}
\end{figure}

The original \TTf{} algorithm is based on consecutive applications of
singular value decompositions (SVDs) on unfoldings of a tensor
\cite{oseledets2011tensor}.
Although known for its efficiency, \TTf{} algorithm requires access to the
full tensor, which is impractical and even impossible for extra-large
tensors.
To address this challenge, the cross interpolation algorithm,
TT-cross, has been developed \cite{oseledets2010tt}.
The idea behind TT-cross is essentially to replace the SVD in the TT
algorithm with an \emph{approximate} version of the skeleton/CUR
decomposition \cite{goreinov1997theory, mahoney2009cur}.
CUR decomposition approximates a matrix by selecting a few of its
columns $\mat{C}$, a few of its rows $\mat{R}$, and a matrix $\mat{U}$
that connects them, as shown in Fig.~\ref{fig:CUR}.
Mathematically, CUR decomposition finds an approximation for a matrix
$\mat{A}$, as $\mat{A}\approx \mat{C}\mat{U}\mat{R}$.
The TT-cross algorithm, utilizes the Maximum Volume Principle
(\emph{maxvol algorithm})
\cite{goreinov2010find,mikhalev2018rectangular} to determine
$\mat{U}$.
The maxvol algorithm chooses a few columns, $\mat{C}$ and rows,
$\mat{R}$, of $\mat{A}$, such that, the intersection matrix
$\mat{U}^{-1}$ has maximum volume \cite{savostyanov2011fast}.

TT-cross interpolation and its versions can be seen as an heuristic
generalization of CUR to tensors
\cite{oseledets2008tucker,sozykin2022ttopt}.
TT-cross utilizes the maximum volume algorithm iteratively, often
beginning with few randomly chosen fibers, to select an optimal number
of specific tensor fibers that capture essential information of the
tensor \cite{savostyanov2014quasioptimality}.
These fibers are used to construct a lower-rank \TTf{}
representation. The naive generalization of CUR is proven to be
expensive, which led to the development of various heuristic
optimization techniques, such as, TT-ALS \cite{holtz2012alternating},
DMRG \cite{oseledets2012solution,savostyanov2011fast}, and AMEN
\cite{dolgov2014alternating}.

\subsection{Quantized Tensor Train format}
\label{sec:QTT}
When the TT-cores of linear operators in TT-format possess a special tensor structures, they can be further compressed. The quantized tensor train (QTT) format is an extension of the tensor train (TT) decomposition, designed to handle high-dimensional data more efficiently by exploiting hidden low-rank structures in exponentially large space. In standard TT decomposition, a high-dimensional tensor is factorized into a product of smaller, lower-dimensional tensors (the TT-cores), reducing the storage complexity and making computations more tractable. However, when dealing with extremely large tensors (e.g., those arising from discretizations with many degrees of freedom), even the TT decomposition may become computationally expensive.

QTT attempts to address this by applying TT decomposition recursively to the modes of the original tensor after reshaping it into a higher-dimensional structure. Specifically, a large tensor is first reshaped into a sequence of smaller tensors using a dyadic (power-of-two) splitting of its dimensions. Each smaller tensor is then factorized using the TT format. This multi-level approach can dramatically reduce storage requirements and computational costs, especially for large structured tensors.  The   QTT format is particularly useful when the underlying data exhibits self-similar patterns or regular structures that can be efficiently captured at multiple scales. For example, if a TT-core has a Toeplitz structure, we can compress
it using \QTT{format}, since Toeplitz matrices have low-rank
\QTT{formats} \cite{kazeev2012low}. This is the case for the CDR transport problem we solve below.
Therefore, to achieve higher compression, we further transform these
TT cores into QTT format, and then solve the mixed TT/QTT
version of the problems with an appropriate TT optimization
technique.



\section{Tensorization of the SEM}
\label{sec:Tensorization}

In this section, we detail the tensorization of the SEM for the space-time CDR Equation. We will begin with deriving the local mass, stiffness, and time derivative matrices associated with the system \ref{Weak:Form}.  From there, the Petrov-Galerkin approximation will be written in terms of a sum of Kronecker products of the local matrices. This will then be used to assemble global matrices for the system.  This ultimately leads to the tensorization of the weak form of the CDR.

\subsection{Local One-dimensional Mass, Stiffness, and Time-derivative Matrices}
We now define the local  mass (time derivative) and stiffness (diffusion) matrices. Since test functions and local basis functions are tensor products of 1D hat functions, these matrices can be described as Kronecker products of local 1D mass and stiffness matrices. From Equations \ref{eq: phi} and \ref{eq: derivative phi}, we see that the local 1D mass and stiffness matrices at index $k$ for the variable $x$ is given by 
\begin{equation}
\bold{M}_{I^x_{k}}= \begin{bmatrix}
\langle \phi_{0,k}^x, \phi_{0,k}^x \rangle & \langle \phi_{1,k}^x, \phi_{0,k}^x \rangle \\
\langle \phi_{0,k}^x, \phi_{1,k}^x \rangle & \langle\phi_{1,k}^x, \phi_{1,k}^x \rangle
\end{bmatrix}, ~
 \bold{S}_{I^x_{k}}= \begin{bmatrix}
\langle\frac{d\phi_{0,k}^x}{dx}, \frac{d\phi_{0,k}^x}{dx}\rangle & \langle\frac{d\phi_{1,k}^x}{dx}, \frac{d\phi_{0,k}^x}{dx}\rangle \\
\langle\frac{d\phi_{0,k}^x}{dx}, \frac{d\phi_{1,k}^x}{dx}\rangle & \langle\frac{d\phi_{1,k}^x}{dx}, \frac{d\phi_{1,k}^x}{dx}\rangle
\end{bmatrix}.
\label{one:stiff:mass}
\end{equation}
The matrices for other space-intervals, $I^t_{k}, I^y_{k},I^z_{k},$ are defined in a similar way. For the local time derivative matrix, $\bold{D}_{I^t_{k}}$, we have the following expression,

\begin{equation}
    \bold{D}_{I^t_{k}}= \begin{bmatrix}
\langle \frac{d\phi_{0,k}^t}{dt}, \phi_{0,k}^t\rangle & \langle \frac{d\phi_{1,k}^t}{dt}, \phi_{0,k}^t\rangle \\
\langle \frac{d\phi_{0,k}^t}{dt}, \phi_{1,k}^t\rangle & \langle\frac{d\phi_{1,k}^t}{dt}, \phi_{1,k}^t\rangle
\end{bmatrix}.
\label{one:time:derivative}
\end{equation}
Importantly, since we are using homogeneous a space-time mesh, the defined mass and stiffness matrices are the same for each four-dimensional hyper-cube and do not depend on the index $k$. 

\subsection{Local Discretization of the Variational Form}
Using these local matrices, we detail the discretization and matricization of each of the bilinear forms $\bold{T_i}$ in Eq.~\eqref{Weak:Form2}. Throughout, we will pick the local basis functions $\phi_{i,l}$ for our test functions $v_h$.



\subsubsection{ Discretization of the time-derivative term, $\bold{T_1}$, on a local four-dimensional hyper-cube}
The discrete approximation of the time-derivative term, on the local four-dimensional hyper-cube, $q_l$, is
\begin{equation}
     \bold{T}_1 := \left \langle  \frac{\partial u^l_h}{\partial t}, \phi_{j} \right \rangle 
     = \sum_{i=0}^{15} u(V_{l,i}) \int_{I^t_{k_t} \times I^x_{k_x} \times I^y_{k_y} \times I^z_{k_z}} \frac{\partial \phi_{i}(t,\bold{x})}{\partial t}  \phi_{j}(t,\bold{x}) \,dt \, dx\, dy\, dz.
      \label{time-der}
 \end{equation}      
 For simplification, we subpress the local element index $l$ from the basis functions and $\bold{T_i}$ unless it is specified.
Using the properties of multiple-integrals on the Cartesian product mesh, and the fact that the space-time local functions are expressed as products of one-dimensional functions, this can be compactly expressed through the one-dimensional mass, stiffness, and time-derivative matrices, using the Kronecker product. When we consecutively use for $v_h$ every one of the $16$ space-time function, $\phi_{j}(t,\bold{x})$, defined on the hyper-cube, $q_l$, Eq.~\eqref{time-der} results in $16$ equations with $16$ unknowns, which in matrix notations can be expressed as:

 \begin{equation}    
  \bold{T}_1=   \Big ( \bold{D}_{I^t_{k_t}} \otimes  \bold{M}_{I^z_{k_z}} \otimes \bold{M}_{I^y_{k_y}} \otimes \bold{M}_{I^x_{k_x}} \Big) \bold{U},
\end{equation}
where $\bold{U}$ is a $16 \times 1$ column-vector containing the unknown $2^4$ components of $u_h$ on $q_l$.
\subsubsection{ Discretization of the diffusion term, $\bold{T_2}$, on four-dimensional hyper-cube, with a constant diffusion coefficient, $\kappa(t,\bold{x})=1$}

In this section, we derive the matrization of $\bold{T_2}$ with $\kappa(t,\bold{x})=1$, so that the term reduces to $\bold{T_2}= \langle  \nabla u_h , \nabla v_h \rangle $. By using \eqref{Lagrange:operator}, and $v_h=\phi_{j}(t,\bold{x})$, we can explicitly write $\bold{T_2}$ locally on $q_l$ as,
   \begin{equation}
       \begin{split}
           \bold{T_2}& := \sum_{i=0}^{15}u(V_{l,i}) \int_{I^t_{k_t} \times I^x_{k_x} \times I^y_{k_y} \times I^z_{k_z}} ( \nabla \phi_{i}(t,\bold{x}) , \nabla \phi_{j}(t,\bold{x}))\,dt \, dx\, dy\, dz  \\
           & = \sum_{i=0}^{15} u(V_{l,i}) \Bigg ( \int_{I^t_{k_t} \times I^x_{k_x} \times I^y_{k_y} \times I^z_{k_z}}  \Big ( \nabla \phi_{i}(t,\bold{x})|{_x} \nabla \phi_{j}(t,\bold{x})|{_x}+\nabla \phi_{i}(t,\bold{x})|{_y} \nabla \phi_{j}(t,\bold{x})|{_y} \\
           & \qquad \qquad +\nabla \phi_{i}(t,\bold{x})|{_z} \nabla \phi_{j}(t,\bold{x})|{_z} \Big)  \,dt \, dx\, dy\, dz  \Bigg).
       \end{split}
       \label{Diff_constant}
   \end{equation}
   This results in $16$ expressions with $16$ unknowns, which in matrix notations becomes:
   \begin{equation}
       \bold{T_2}:= \big (\bold{M}_{I^t_{k_t}} \otimes \bold{S}_{I^z_{k_z}} \otimes \bold{M}_{I^y_{k_y}} \otimes \bold{M}_{I^x_{k_x}}  +\bold{M}_{I^t_{k_t}} \otimes \bold{M}_{I^z_{k_z}} \otimes \bold{S}_{I^y_{k_y}} \otimes \bold{M}_{I^x_{k_x}}  + \bold{M}_{I^t_{k_t}} \otimes \bold{M}_{I^z_{k_z}} \otimes \bold{M}_{I^y_{k_y}} \otimes \bold{S}_{I^x_{k_x}} \big) \bold{U}.
       \label{Space-time_diff_cons}
   \end{equation}
\subsubsection{ Discretization of the diffusion term, $\bold{T_2}$, on four-dimensional hyper-cube, with an non-constant diffusion function}
For variable diffusion coefficient, $\kappa(t, \bold{x})$, we employ the Lagrange interpolation operator defined in \eqref{Lagrange:operator}. Assuming that $\kappa(t, \bold{x})$ is well-defined at the nodes of $q_l$, and once again using $v_h = \phi_j$, we derive the matrization of $\bold{T_2}$ as 
\begin{equation}
       \begin{split}
           \bold{T_2} & = \langle L_h^l (\kappa(t,\bold{x})) \nabla u_h^l, \nabla v_h \rangle  =  \sum_{i=0}^{15}u(V_{l,i})  \int_{q_l} L_h^l (\kappa(t,\bold{x})) ( \nabla \varphi_{i}(t,\bold{x}) , \nabla \varphi_{j}(t,\bold{x}) )dtdxdydz \\
           &= \sum_{i,m=0}^{15} \kappa(V_{l,m})  u(V_{l,i}) \int_{I^t_{k_t} \times I^x_{k_x} \times I^y_{k_y} \times I^z_{k_z}}   \phi_n(t,\bold{x}) \Big (  \nabla \phi_{i}(t,\bold{x})|{_x} \nabla \phi_{j}(t,\bold{x})|{_x} \\
           & \qquad +\nabla \phi_{i}(t,\bold{x})|{_y} \nabla \phi_{j}(t,\bold{x})|{_y}  +  \nabla \phi_{i}(t,\bold{x})|{_z} \nabla \phi_{j}(t,\bold{x})|{_z} \Big)  \,dt \, dx\, dy\, dz ,
       \end{split}
       \label{Diff_var}
   \end{equation}
where $\kappa(V_{l,m})$ are the values of the diffusion function $\kappa(t, \bold{x})$ at the $2^4$ nodes of $q_l$. Here, in contrast to \eqref{Diff_constant}, each one-dimensional integration is on the product of one basis function and two gradients of basis functions. To formulate this in a way that is structurally similar to \eqref{Space-time_diff_cons}, we introduce for each dimension following modified mass and stiffness matrices
\begin{equation}
\bold{M}_{I^x_k}^{p}= \begin{bmatrix}
\langle \phi_{p,k}^x \phi_{0,k}^x, \phi_{0,k}^x \rangle & \langle \phi_{p,k}^x \phi_{1,k}^x, \phi_{0,k}^x \rangle  \\
\langle  \phi_{p,k}^x \phi_{0,k}^x, \phi_{1,k}^x \rangle & \langle \phi_{p,k}^x \phi_{1,k}^x, \phi_{1,k}^x\rangle
\end{bmatrix}
;\bold{S}_{I^x_k}^{p}= \begin{bmatrix}
\langle  \phi_{p,k}^x \frac{d\phi_{0,k}^x}{dx}, \frac{d\phi_{0,k}^x}{dx} \rangle &  \langle \phi_{p,k}^x \frac{d\phi_{1,k}^x}{dx}, \frac{d\phi_{0,k}^x}{dx}\rangle  \\
\langle \phi_{p,k}^x \frac{d\phi_{0,k}^x}{dx}, \frac{d\phi_{1,k}^x}{dx} \rangle  & \langle \phi_{p,k}^x \frac{d\phi_{1,k}^x}{dx}, \frac{d\phi_{1,k}^x}{dx}\rangle 
\end{bmatrix},
\label{Loc_mat_var}
\end{equation}
where the index $p$ assumes two values: $0$, or $1$.
For $v_h= \phi_j$, Equation \eqref{Diff_var} reduces to a matrix-vector product, where the associated matrix can be constructed by using one-dimensional matrices introduced in \eqref{Loc_mat_var}. Using again the properties of multiple-integrals on the Cartesian product mesh, and the fact that the
space-time local functions are products of one-dimensional functions, explicitly, we have, 
\begin{equation}
\begin{split}
    \bold{T_2}&:= \Bigg ( \sum_{m_1,m_2,m_3,m_4=0}^{1} \kappa(V_{l,m}) (  \bold{M}_{I^t_{k_t}}^{{m_4}} \otimes \bold{S}_{I^z_{k_z}}^{{m_3}} \otimes \bold{M}_{I^y_{k_y}}^{{m_2}} \otimes \bold{M}_{I^x_{k_x}}^{{m_1}}   
 \\
 &+\bold{M}_{I^t_{k_t}}^{{m_4}} \otimes \bold{M}_{I^z_{k_z}}^{{m_3}} \otimes \bold{S}_{I^y_{k_y}}^{{m_2}} \otimes \bold{M}_{I^x_{k_x}}^{{m_1}}   + \bold{M}_{I^t_{k_t}}^{{m_4}} \otimes \bold{M}_{I^z_{k_z}}^{{m_3}} \otimes \bold{M}_{I^y_{k_y}}^{{m_2}} \otimes \bold{S}_{I^x_{k_x}}^{{m_1}} )\Bigg ) \bold{U}.
 \end{split}
 \label{T2:mat}
\end{equation}

\subsubsection{ Discretization of the convection term, $\bold{T_3}$, on four-dimensional hyper-cube, with an non-constant convection function}
 In this section, we derive the matrization of the convection term with variable coefficients. Following the same approach as in \eqref{Diff_var}  
using \eqref{Dis_uh}, and $v_h=\phi_{j}(t,\bold{x})$, $\bold{T_3}$ becomes
\begin{equation}
    \begin{split}
        \bold{T_3}&=\langle L_h^l(\bold{b}(t,\bold{x})) \cdot \nabla u_h^l, v_h \rangle \\
        &= \sum_{m=0}^{15}b_1(V_{l,m}) \sum_{i=0}^{15} u(V_{l,i})  \int_{q_l}   \phi_m(t,\bold{x})   \nabla \phi_{i}(t,\bold{x})|{_x}  \phi_{j}(t,\bold{x}) dtdxdydz\\
           & \quad +\sum_{m=0}^{15}b_2(V_{l,m})  \sum_{i=0}^{15} u(V_{l,i})  \int_{q_l}   \phi_m(t,\bold{x})   \nabla \phi_{i}(t,\bold{x})|{_y}  \phi_{j}(t,\bold{x}) dtdxdydz \\
           & \quad + \sum_{m=0}^{15}b_3(V_{l,m}) \sum_{i=0}^{15} u(V_{l,i})  \int_{q_l}   \phi_m(t,\bold{x})   \nabla \phi_{i}(t,\bold{x})|{_z}  \phi_{j}(t,\bold{x}) \,dt \, dx\, dy\, dz.
    \end{split}
    \label{T_3:Int}
\end{equation}
To derive the matrix structure of \eqref{T_3:Int}, we introduce the following one-dimensional local space-derivative matrix for each dimension. On the interval $I^x_l$, we define
\begin{equation}
 \bold{D}_{I^x_k}^{p}= \begin{bmatrix}
\langle \phi_{p,k}^x \frac{d\phi_{0,k}^x}{dx}, \phi_{0,k}^x \rangle &  \langle \phi_{p,k}^x \frac{d\phi_{1,k}^x}{dx}, \phi_{0,k}^x \rangle \\
\langle \phi_{p,k}^x \frac{d\phi_{0,k}^x}{dx}, \phi_{1,k}^x \rangle & \langle \phi_{p,k}^x \frac{d\phi_{1,k}^x}{dx}, \phi_{1,k}^x \rangle
\end{bmatrix},
\label{T_3:local:x}
\end{equation}
where again $p\in \{0,1\}$.  By using \eqref{T_3:local:x}, and \eqref{Loc_mat_var}, we find that
\begin{equation}
\begin{split}
    \bold{T_3} 
&=\Bigg (\sum_{m_1,m_2,m_3,m_4=0}^1  b_1(V_{l,m})  (\bold{M}_{I^t_{k_t}}^{{m_4}}\otimes \bold{M}_{I^z_{k_z}}^{{m_3}}\otimes \bold{M}_{I^y_{k_y}}^{{m_2}} \otimes \bold{D}_{I^x_{k_x}}^{{m_1}} )
\\
& \qquad \qquad + b_2(V_{l,m})  (\bold{M}_{I^t_{k_t}}^{{m_4}} \otimes \bold{M}_{I^z_{k_z}}^{{m_3}} \otimes \bold{D}_{I^y_{k_y}}^{{m_2}} \otimes \bold{M}_{I^x_{k_x}}^{{m_1}} )   
\\
& \qquad \qquad + b_3(V_{l,m})  (\bold{M}_{I^t_{k_t}}^{{m_4}} \otimes \bold{D}_{I^z_{k_z}}^{{m_3}} \otimes \bold{M}_{I^y_{k_y}}^{{m_2}} \otimes \bold{M}_{I^x_{k_x}}^{{m_1}} )   \Bigg ) \bold{U}.
    \end{split}
\end{equation}

\subsubsection{ Discretization of the reaction term, $\bold{T_4}$, on four-dimensional hyper-cube}
In this section, we build the matrization of the reaction term, $\bold{T_4}$, with variable reaction coefficient. We highlight that for a constant reaction coefficient and $v_h = \phi_j$, the matrix associated with $\bold{T_4}$ reduces to the mass matrix.
In this case, it follows from \eqref{one:stiff:mass} that
\begin{equation}
\begin{split}
\bold{T_4}  = \Big( \bold{M}_{I^t_{k_t}} \otimes \bold{M}_{I^x_{k_x}} \otimes \bold{M}_{I^y_{k_y}} \otimes \bold{M}_{I^x_{k_x}}\Big ) \bold{U}.
\end{split}
\end{equation}
For variable reaction coefficient, we discretize the reaction term using Lagrange interpolation operator \eqref{Lagrange:operator}. Again taking $v_h=\phi_{j}(t,\bold{x})$, $\bold{T_4}$ becomes
\begin{equation*}
    \bold{T_4}=\langle L_h^l(c(t,\bold{x})) u_h^l, v_h \rangle 
     = \sum_{m,i=0}^{15} c(V_{l,m}) u(V_{l,i})  \int_{q_l}   \phi_m(t,\bold{x})    \phi_{i}(t,\bold{x})  \phi_{j}(t,\bold{x}) ~dtdxdydz,
\end{equation*}
Using \eqref{Loc_mat_var}, we derive
\begin{equation}
    \bold{T_4}= \Bigg (\sum_{m_1,m_2,m_3,m_4=0}^1  c(V_{l,m}) (\bold{M}_{I^t_{k_t}}^{{m_4}} \otimes \bold{M}_{I^z_{k_z}}^{{m_3}} \otimes \bold{M}_{I^y_{k_y}}^{{m_2}} \otimes \bold{M}_{I^x_{k_x}}^{{m_1}}) \Bigg )\bold{U}.
\end{equation}
\subsubsection{ Discretization of the loading term, $\bold{T_5}$, on four-dimensional hyper-cube}
Finally, using Lagrange interpolation operator \eqref{Lagrange:operator}, we rewrite the load vector into matrix vector product. For $v_h=\phi_{j}(t,\bold{x})$, $\bold{T_5}$ reduces to
\begin{equation*}
    \bold{T_5}=\langle L_h^l(f(t,\bold{x})) , v_h \rangle 
    = \sum_{m=0}^{15} f(V_{l,m}) \int_{q_l}   \varphi_m(t,\bold{x})   \varphi_{j}(t,\bold{x}) \,dt \, dx\, dy\, dz.
\end{equation*}
The final matrix vector product can be written as 
\begin{equation}
    \bold{T_5}:= (\bold{M}_{I^t_{k_t}} \otimes \bold{M}_{I^z_{k_z}} \otimes \bold{M}_{I^y_{k_y}} \otimes \bold{M}_{I^x_{k_x}} ) \bold{F}.
\end{equation}
Here, the load vector is defined as 
\[
(\bold{F})_m:=f(V_{l,m}) ,
\]
where $m$ (as always) denotes the $16$ local nodes of $q_l$.

\section{ Assembly of Global Matrices}
\label{sec:aseemble}
\subsection{ Assembly for Terms with Constant Coefficients}
The local matrix structures of $\bold{T_1}, \bold{T_2}, \bold{T_3}, \bold{T_4}$, and $\bold{T_5}$ can easily be extended to a similar global structures when the diffusion, convection, and reaction coefficients are constants. To achieve this we need to construct the global: mass, stiffness, and time-derivative matrices, using the local ones for each variable.

First, we construct from  \eqref{one:time:derivative} the following block matrix, $\bold{D}_{t}$, for the temporal variable.  Recall that the time interval $I_T$ is decomposed using uniform intervals $I_T = \cup_{k_t=1}^N I_{k_t}^t$.  We define
\begin{equation}
\bold{D}_{t} = \begin{bmatrix}
  \bold{D}_{I^t_1} & 0 &0 &0   &\ldots  & 0\\
  0 & \bold{D}_{I^t_2} &0 &0 &\ldots &0 \\
  \vdots & \vdots & 
  \vdots & \vdots &\vdots & \vdots \\
  0 &0 &0 &\ldots &\bold{D}_{I^t_{N-1}} &0 \\
  0&0 &0 &0 &0 &\bold{D}_{I^t_{N}}
\end{bmatrix}_{2N\times 2N},
\label{B:Omega}
\end{equation}
where $\bold{D}_{I_{k_t}^t}$ is as defined in Equation \ref{one:time:derivative}.
 To achieve the proper global assembly of the 1D time derivative matrix, we need to take into account the common boundary points of each interval. This can be accomplished with the help of the binary matrix, $\bold{B}$
\begin{equation}
\bold{B}= \begin{bmatrix}
  1 & 0 &0 &0   &0 &\ldots  & 0\\
  0 & 1 &1 &0 &0 & \ldots  &0 \\
  0 &0 &0 &1 &1 &  \ldots  &0\\
  \vdots  & \vdots & \vdots &\vdots & \vdots & \vdots \\
  0&0 &0 &0 &0 &\ldots  &1
\end{bmatrix}_{(N+1)\times 2N}.
\label{B:L}
\end{equation}
Then, the global $[(N+1) \times (N+1)]$ time-derivative matrix, $\bold{D}_{{I_T}}$, can be determined as
\begin{equation}
    \bold{D}_{I_T}= \bold{B}\bold{D}_{t} \bold{B}^T.
\end{equation}
In a similar fashion, using $\bold{B}$, we can build the global one-dimensional mass matrices $\bold{M}_{\Omega_Z}, \bold{M}_{\Omega_Y}$, and $\bold{M}_{\Omega_X} $. Finally, we obtain for $\bold{T}_1^{\text{g}}$, 
\begin{equation}
\label{T1:global}
  \bold{T}_1^{\text{g}}:=   \Big ( \bold{D}_{{I_T}} \otimes  \bold{M}_{\Omega_Z} \otimes \bold{M}_{\Omega_Y} \otimes \bold{M}_{\Omega_X} \Big) \bold{U}^{\text{g}},
\end{equation}
here $\bold{U}^{\text{g}}$ is a global column solution vector of size $[(N+1)^4 \times 1] $. Similarly, we build the term $\bold{T}_5^{\text{g}}$, 

\begin{equation}
\label{T5:global}
    \bold{T}_5^{\text{g}}:= (\bold{M}_{I_T} \otimes \bold{M}_{\Omega_Z} \otimes \bold{M}_{\Omega_Y} \otimes \bold{M}_{\Omega_X} ) \bold{F}^{\text{g}},
\end{equation}
where $\bold{F}^{\text{g}}$ is the column load vector of size $ [(N+1)^4 \times 1]$. 
\subsection{ Assembly for Terms with Variable Coefficients}

\subsubsection{Assembly of $\mathbf{T_2}$, when diffusion coefficient depends on space-time variables, but allows a separation of variables}
First, let us assume that, 
\begin{equation}
\kappa(t,z,y,x)= \kappa^t(t) \kappa^z(z) \kappa^y(y) \kappa^x(x),
\end{equation}
$\kappa^t(t)$, $\kappa^z(z)$, $\kappa^y(y)$, and $\kappa^x(x)$ are separation of variables of the function $\kappa(t,z,y,x)$, then we have locally (\ref{T2:mat}),

\begin{equation}
\begin{split}
    \bold{T}_2& = \bigg (\sum_{m_1,m_2,m_3,m_4=0}^{1} \kappa^t(t_{m_4})\kappa^z({z_{m_3}})\kappa^y({y_{m_2}})\kappa^x({x_{m_1}})(  \bold{M}_{I^t_{k_t}}^{{m_4}} \otimes \bold{S}_{I^z_{k_z}}^{{m_3}} \otimes \bold{M}_{I^y_{k_y}}^{{m_2}} \otimes \bold{M}_{I^x_{k_x}}^{{m_1}}   
 \\
 &+\bold{M}_{I^t_{k_t}}^{{m_4}} \otimes \bold{M}_{I^z_{k_z}}^{{m_3}} \otimes \bold{S}_{I^y_{k_y}}^{{m_2}} \otimes \bold{M}_{I^x_{k_x}}^{{m_1}}   + \bold{M}_{I^t_{k_t}}^{{m_4}} \otimes \bold{M}_{I^z_{k_z}}^{{m_3}} \otimes \bold{M}_{I^y_{k_y}}^{{m_2}} \otimes \bold{S}_{I^x_{k_x}}^{{m_1}} )\Bigg ) \bold{U}.
 \end{split} 
\end{equation}   
To assemble the global matrix in $x$-dimension, let us define the following pairs of local one-dimensional coefficient matrices. For each interval $I^x_{k-1} = [x_{k-1}, x_k]$, we have:
\begin{equation}
\begin{split}
\boldsymbol{\kappa}^{x}_{{k-1}}= \begin{bmatrix}
\kappa^x(x_{k-1}) &  0 \\
0 & \bold{\kappa}^x(x_{k})
\end{bmatrix}_{2 \times 2}.
\end{split}
\end{equation}
and similarly for $y,z$, and $t$ dimensions. By evaluating $\kappa^x(x)$ at all points of $\Omega_X=\cup_{k_x=1}^{N} I^x_{k_x}$, we define the following pair of global block coefficient matrices in $x$-dimension,
\begin{equation}
\bold{C}^{x}_{{0}}= \begin{bmatrix}
  \boldsymbol{\kappa}^{x}_{{0}} & 0 &0 &0   &\ldots  & 0\\
  0 & \boldsymbol{\kappa}^{x}_{{1}} &0 &0 &\ldots &0 \\
  \vdots & \vdots & 
  \vdots & \vdots &\vdots & \vdots \\
  0 &0 &0 &\ldots &\boldsymbol{\kappa}^{x}_{{{N-2}}} &0 \\
  0&0 &0 &0 &0 &\boldsymbol{\kappa}^{x}_{{{N-1}}}
\end{bmatrix}_{2N \times 2N},
\label{C:Omega:L-1}
\end{equation}

\begin{equation}
\bold{C}^{x}_{{1}}= \begin{bmatrix}
  \boldsymbol{\kappa}^{x}_{{1}} & 0 &0 &0   &\ldots  & 0\\
  0 & \boldsymbol{\kappa}^{x}_{{2}} &0 &0 &\ldots &0 \\
  \vdots & \vdots & 
  \vdots & \vdots &\vdots & \vdots \\
  0 &0 &0 &\ldots &\boldsymbol{\kappa}^{x}_{{{N-1}}} &0 \\
  0&0 &0 &0 &0 &\boldsymbol{\kappa}^{x}_{{{N}}}
\end{bmatrix}_{2N \times 2N},
\label{C:Omega:L}
\end{equation}
Using \eqref{Loc_mat_var}, we construct the  global pair of stiffness block matrices as follows,

\begin{equation}
\bold{S}^x_{0}= \begin{bmatrix}
  \bold{S}^{0}_{I^{x}_{1}} & 0 &0 &0   &\ldots  & 0\\
  0 & \bold{S}^{0}_{I^x_{2}} &0 &0 &\ldots &0 \\
  \vdots & \vdots & 
  \vdots & \vdots &\vdots & \vdots \\
  0 &0 &0 &\ldots &\bold{S}^{{0}}_{I^x_{N-1}} &0 \\
  0&0 &0 &0 &0 &\bold{S}^{{0}}_{I^x_{N}}
\end{bmatrix}_{2N \times 2N},
\label{S:Omega:L-1}
\end{equation}

\begin{equation}
\bold{S}^x_1= \begin{bmatrix}
  \bold{S}^{{1}}_{I^x_1} & 0 &0 &0   &\ldots  & 0\\
  0 & \bold{S}^{{1}}_{I^x_2} &0 &0 &\ldots &0 \\
  \vdots & \vdots & 
  \vdots & \vdots &\vdots & \vdots \\
  0 &0 &0 &\ldots &\bold{S}^{x_{N-1}}_{I^x_{N-1}} &0 \\
  0&0 &0 &0 &0 &\bold{S}^{N}_{I^{x}_{N}}
\end{bmatrix}_{2N \times 2N},
\label{S:Omega:L}
\end{equation}
and similarly for $\bold{M}^x_0$ and $\bold{M}^x_{1}$. By using \eqref{B:L},~\eqref{C:Omega:L-1}, and \eqref{S:Omega:L-1} we construct the global one-dimensional matrices,
\begin{equation}
    \bold{S}^{\text{g}}_{x_{0}}:= \bold{B}\bold{S}^x_{0}\bold{C}^{x}_{{0}} \bold{B}^T, \text{and }\bold{S}^{\text{g}}_{x_{1}}:= \bold{B}\bold{S}^x_{1}\bold{C}^{x}_{{1}} \bold{B}^T,
    \label{Glb:S^g}
\end{equation}
and similarly, the global stiffness matrices $\bold{S}^{\text{g}}_{y_{0}}$, $\bold{S}^{\text{g}}_{z_{0}}$, $\bold{S}^{\text{g}}_{y_{1}}$, and $\bold{S}^{\text{g}}_{z_{1}}$. Next, for the $x$-dimensional global mass matrix we have,
\begin{equation}
    \bold{M}^{\text{g}}_{x_{0}}:= \bold{B}\bold{M}^x_{0}\bold{C}^{x}_{{0}} \bold{B}^T,\text{and } \bold{M}^{\text{g}}_{x_{1}}:= \bold{B}\bold{M}^x_{1}\bold{C}^{x}_{{1}} \bold{B}^T,
\end{equation}
and similarly for the remaining dimensions. By using the above matrices, we construct the global space-time matrization for diffusion term as,
\begin{equation}
\begin{split}
    \bold{T}_2^{\text{g}}&:= \Bigg ( \sum_{m_1,m_2,m_3,m_4=0}^{1}   [\bold{M}^{\text{g}}_{t_{m_4}} \otimes \bold{S}^{\text{g}}_{z_{m_3}} \otimes \bold{M}^{\text{g}}_{y_{m_2}} \otimes \bold{M}^{\text{g}}_{x_{m_1}}   
 \\
 &+\bold{M}^{\text{g}}_{t_{m_4}} \otimes \bold{M}^{\text{g}}_{z_{m_3}} \otimes \bold{S}^{\text{g}}_{y_{m_2}} \otimes \bold{M}^{\text{g}}_{x_{m_1}}   + \bold{M}^{\text{g}}_{t_{m_4}} \otimes \bold{M}^{\text{g}}_{z_{m_3}} \otimes \bold{M}^{\text{g}}_{y_{m_2}} \otimes \bold{S}^{\text{g}}_{x_{m_1}}] \Bigg ) \bold{U^{\text{g}}}.
 \end{split}
\label{Diff:coef:Omega}
\end{equation}
%

\subsection{ Assembly for convection term $\bold{T_3}$, on four-dimensional hyper-cube}

Proceeding similarly as \eqref{Diff:coef:Omega}, we derive the low-rank structure of the convective term, and reaction terms as follows

\begin{equation}
\begin{split}
    \bold{T}_3^{\text{g}} 
&=\Bigg (\sum_{m_1,m_2,m_3,m_4=0}^{1}    [\bold{M}^{\text{g}}_{t_{m_4}}\otimes \bold{M}^{\text{g}}_{z_{m_3}}\otimes \bold{M}^{\text{g}}_{y_{m_2}} \otimes \bold{D}^{\text{g}}_{x_{m_1}} 
\\
& \qquad \qquad +   \bold{M}^{\text{g}}_{t_{m_4}}\otimes \bold{M}^{\text{g}}_{z_{m_3}}\otimes \bold{D}^{\text{g}}_{y_{m_2}} \otimes \bold{M}^{\text{g}}_{x_{m_1}}  
\\
& \qquad \qquad +  \bold{M}^{\text{g}}_{t_{m_4}}\otimes \bold{D}^{\text{g}}_{z_{m_3}}\otimes \bold{M}^{\text{g}}_{y_{m_2}} \otimes \bold{M}^{\text{g}}_{x_{m_1}}]  \Bigg ) \bold{U}^{\text{g}},
    \end{split}
    \label{Convec:coef:Omega}
\end{equation}
where
\begin{equation}
    \bold{D}^{\text{g}}_{x_{0}}:= \bold{B}\bold{D}^x_{0}\bold{C}^{x}_{{0}} \bold{B}^T \text{and } \bold{D}^{\text{g}}_{x_{1}}:= \bold{B}\bold{D}^x_{1}\bold{C}^{x}_{{1}} \bold{B}^T.
    \label{Glb:M^g}
\end{equation}
The matrices $\bold{D}^{\text{g}}_{y_{0}}$, $\bold{D}^{\text{g}}_{y_{1}}$, $\bold{D}^{\text{g}}_{z_{0}}$ and $\bold{D}^{\text{g}}_{z_{1}}$ are constructed by using \eqref{T_3:local:x}, \eqref{B:Omega} and the similar steps used to build \eqref{Glb:S^g}.

\subsection{ Assembly for reaction term $\bold{T_4}$ on four-dimensional hyper-cube}
By using $\bold{M}^{\text{g}}_{t_{m_4}}$, $\bold{M}^{\text{g}}_{z_{m_3}}$, $\bold{M}^{\text{g}}_{y_{m_2}}$, and $\bold{M}^{\text{g}}_{x_{m_1}}$, we can construct the matrization corresponding to reaction term $\bold{T_4}$ as
\begin{equation}
\label{T4:global}
    \bold{T}^{\text{g}}_4:= \sum_{m_1,m_2,m_3,m_4=0}^{1}   (\bold{M}^{\text{g}}_{t_{m_4}} \otimes \bold{M}^{\text{g}}_{z_{m_3}} \otimes \bold{M}^{\text{g}}_{y_{m_2}} \otimes \bold{M}^{\text{g}}_{x_{m_1}}) \bold{U}^{\text{g}}.
\end{equation}


\subsection{The Global System} 
Throughout, we have assumed that $\kappa$ was a rank one function in order to derive  expressions for $\bold{T}^{\text{g}}_1,\bold{T}^{\text{g}}_2, \bold{T}^{\text{g}}_3, \bold{T}^{\text{g}}_4$, and $\bold{T}^{\text{g}}_5$.
In the general case, we utilize the cross-interpolation technique~\cite{oseledets2010approximation} to approximate the TT format of general functions, which can be interpreted as an approximation in separated-variable form. 
Hence, the techniques described above can be used to construct the global matrices.  Taking into account all expressions $\bold{T}^{\text{g}}_1,\bold{T}^{\text{g}}_2, \bold{T}^{\text{g}}_3, \bold{T}^{\text{g}}_4$, and $\bold{T}^{\text{g}}_5$, we result in the following system of equation:
\begin{equation}
    (\bold{T}^{\text{g}}_1+\bold{T}^{\text{g}}_2+\bold{T}^{\text{g}}_3+\bold{T}^{\text{g}}_4) := \bold{A}^{\text{g}}\bold{U}^{\text{g}} =\bold{T}^{\text{g}}_5.
    \label{eqn:full linear system}
\end{equation}
Finally, to deal with boundary conditions, we reformulation the system in Eqn.~\eqref{eqn:full linear system} to only include interior nodes
\begin{equation}
    \bold{A}^{g,int}\bold{U}^{g,int} =\bold{T}^{g,int}_5 - \bold{F}^{bd}.
    \label{eqn: reduced linear system}
\end{equation}
where $\bold{F}^{bd}$ is the boundary term incorporating the boundary and initial conditions. More details about this formulation is described in~\cite{adak2024tensor}.

\if 1<0
\subsubsection{Non-homogeneous Initial and Boundary Conditions} 
\label{Bd_cond}


So far, we have not incorporated the boundary conditions (BC) and the
initial condition (IC), given in \eqref{Lin.Sys}, into the
linear system $\bold{A}^g\bold{U}^g = \bold{T}^g_5$ \eqref{cheb_lin_syst}.
In the space-time method, we consider the IC equivalent to the BC.
The nodes of Cartesian SEM grid are split into two parts: $(i)$ BC/IC
nodes, and $(ii)$ interior nodes, where the solution is unknown.
Let $i^{\text{Bd}}$, and $i^{\text{Int}}$ be the set of multi-indices,
the BC/IC and interior nodes, respectively.

Boundary and initial conditions are imposed by explicitly enforcing
the BC/IC nodes to be equal to the BC, $g(t,\bold{x})$, or to the IC, $h(0,\bold{x})$,
and then reducing the linear system for all nodes into a smaller
system only for the interior nodes,
\begin{equation}
  \label{eqn:reduced_system}
  \bold{A}^{\text{Int}} \bold{U}^{\text{Int}}=\bold{F}^{\text{Int}} -\bold{F}^{\text{Bd}}.
\end{equation}

To make it clear we consider below a simple example with $N = 4$
collocation points in two dimensions $(t,x)$.
The set of Chebyshev nodes can be denoted by multi-indices as,
$\bold{U}:=\{u_1,u_2,\ldots,u_{16} \} $ (Figure
\ref{fig:ST-2d-grid}). After imposing the BC/IC conditions, the linear
system with BC/IC nodes, $\bold{A}\bold{U} = \bold{F}$, is

\begin{equation}
  \label{Full_system}
  \begin{split}
    &u_1 = g (t_{0},{x}_{0})\\
    &\quad \vdots \\
    &u_4 = g (t_{3},{x}_{0})\\
    &u_5 = h(x_1), \\
    &A_{6,1} u_1 + A_{6,2} u_2+ \ldots +A_{6,15}u_{15} + A_{6,16}u_{16} =F_6\\
    &\quad \vdots \\
    &A_{8,1} u_1 + A_{8,2} u_2+ \ldots +A_{8,15}u_{15} + A_{8,16}u_{16} =F_8\\
    &u_9 = h(x_{2})\\
    &A_{10,1} u_1 + A_{10,2} u_2+ \ldots +A_{10,15}u_{15} + A_{10,16}u_{16} =F_{10}\\
    &\quad \vdots \\
    &A_{12,1} u_1 + A_{12,2} u_2+ \ldots +A_{12,15}u_{15} + A_{12,16}u_{16} =F_{12}\\
    & \quad \vdots \\
    &u_{13} = g (t_{0},{x}_{3})\\
    &\quad \vdots \\
    &u_{16} = g (t_{3},{x}_{3}).\\
  \end{split}
\end{equation}

The unknown values $(u_6,u_7,u_8,u_{10},u_{11},u_{12})$ associated to
six interior nodes, satisfy the reduced linear system whose equations
read as
\begin{multline}
  \qquad
  A_{l,6} u_6+A_{l,7}u_7+A_{l,8}u_{8} + A_{l,10}u_{10} +A_{l,11}u_{11} +A_{l,12}u_{12}\\
  =F_l-A_{l,1}g (t_0,x_0)-A_{l,2}g (t_1,x_0)-\ldots,
  \qquad
  \label{eqn:sub_system_BC}
\end{multline}
with $l \in \{6,7,8,10,11,12\}$, where we transfer the values of BC/IC
nodes to the right-hand side.
\fi

\section{Tensorization of the Weak-form of CDR}
\label{sec:Tensor-Networks}
\subsection{Transformation of the CDR Discretization into TT and QTT format}
\label{sec:CDR in TT format}

To solve the CDR equation with coefficients: $\kappa(t,\bold{x})$,
$\bold{b}(t,\bold{x})$, $c(t,\bold{x})$, boundary conditions, initial
conditions, and loading functions, that do not have separation of variables, we use TT-cross to build the TT format,
directly from them.
The space-time discretization process leads to the formulation of a linear system for all interior nodes, as outlined in Eq.~\eqref{eqn: reduced linear system}. For simplicity, further we will omit part of the indices and refer to this equation as: $\bold{A}\bold{U}=\bold{T}_5 -\bold{F}^{\text{bd}}$, where $\bold{A}$ represents the operator matrix, $\bold{U}$ is the solution vector, $\bold{T}_5$ is the loading term, and $\bold{F}^{\text{bd}}$ accounts for the boundary and initial conditions.
In the following section, we will provide a detailed explanation of how to construct the TT and QTT formats for each component of this linear system, utilizing the following three steps: 
\begin{equation}
\begin{split}
  &\bold{A}\bold{U}=\bold{T}_5 -\bold{F}^{\text{bd}}\\ 
  &\ten{A}^{TT}\ten{U}^{TT} = \ten{T}_5^{TT} - \ten{F}^{\text{bd},TT}\\
  &\ten{A}^{QTT}\ten{U}^{QTT} = \ten{T}_5^{QTT} - \ten{F}^{\text{bd},QTT},
\end{split}
  \label{eqn:TT/QTT linear system}
\end{equation}
where: $\ten{A}^{TT}\ten{U}^{TT} := \ten{T}_1^{TT} + \ten{T}_2^{TT} + \ten{T}_3^{TT}+ \ten{T}_4^{TT}$, and 
$\ten{A}^{QTT}\ten{U}^{QTT} := \ten{T}_1^{QTT} + \ten{T}_2^{QTT} + \ten{T}_3^{QTT}+ \ten{T}_4^{QTT}$, are using $\mat{A}, \bold{U}$, and the loading terms - $\bold{T}_5, \bold{F}^{\text{bd}}$ expressed in $TT$ and $QTT$ tensor formats.

\subsection{Construction of TT Format for Linear System Components}
Given that the operators in their matrix representations exhibit a Kronecker product structure, derived in Sec \ref{sec:SEM model}, their TT format can be constructed by employing component matrices as TT cores, see Eq. \eqref{eqn:TT_matrices}. 
To construct the TT format of the operators acting on the \emph{interior}
nodes, we first need to define some sets of indices:
\begin{equation*}
  \begin{split}
    & \ten{I}_t = 2:N+1 \ \text{index set for time variable},\\
    & \ten{I}_s = 2:N \ \text{index set for space variable}.
  \end{split}
\end{equation*}

\PGRAPH{$\bullet$ TT format of $\mathbf{T}^{g,int}_1$, $\ten{T}^{TT}_1$}: From the formulation in Eqn.~\eqref{T1:global}, the global temporal operator in TT-matrix format acting only on the \emph{interior} nodes is constructed as:
\begin{equation}
\label{T1-TT:global}
  \ten{T}_1^{TT}:=   \Big ( \bold{D}_{{I_T}}(\ten{I}_t,\ten{I}_t) \tenprod  \bold{M}_{\Omega_Z}(\ten{I}_s,\ten{I}_s) \tenprod \bold{M}_{\Omega_Y}(\ten{I}_s,\ten{I}_s) \tenprod \bold{M}_{\Omega_X}(\ten{I}_s,\ten{I}_s) \Big) \ten{U}^{TT},
\end{equation}

where the tensor product operator $\tenprod$ is defined in~\ref{appdef}.

\PGRAPH{$\bullet$ TT format of $\mathbf{T}^{g,int}_2$, $\ten{T}^{TT}_2$}: From the formulation in Eqn.~\eqref{Diff:coef:Omega}, the diffusion  operator in TT-matrix format is constructed as:
\begin{equation}
\begin{split}
    \ten{T}_2^{TT}:= \Bigg ( \sum_{m_1,m_2,m_3,m_4=0}^{1}   \big[&\bold{M}^{\text{g}}_{t_{m_4}}(\ten{I}_t,\ten{I}_t) \tenprod \bold{S}^{\text{g}}_{z_{m_3}}(\ten{I}_s,\ten{I}_s) \tenprod \bold{M}^{\text{g}}_{y_{m_2}}(\ten{I}_s,\ten{I}_s) \tenprod \bold{M}^{\text{g}}_{x_{m_1}}(\ten{I}_s,\ten{I}_s)
 \\
 +&\bold{M}^{\text{g}}_{t_{m_4}}(\ten{I}_t,\ten{I}_t) \tenprod \bold{M}^{\text{g}}_{z_{m_3}}(\ten{I}_s,\ten{I}_s) \tenprod \bold{S}^{\text{g}}_{y_{m_2}}(\ten{I}_s,\ten{I}_s) \tenprod \bold{M}^{\text{g}}_{x_{m_1}}(\ten{I}_s,\ten{I}_s)\\
 +&\bold{M}^{\text{g}}_{t_{m_4}}(\ten{I}_t,\ten{I}_t) \tenprod \bold{M}^{\text{g}}_{z_{m_3}}(\ten{I}_s,\ten{I}_s) \tenprod \bold{M}^{\text{g}}_{y_{m_2}}(\ten{I}_s,\ten{I}_s) \tenprod \bold{S}^{\text{g}}_{x_{m_1}}(\ten{I}_s,\ten{I}_s)\big] \Bigg ) \ten{U}^{TT}.
 \end{split}
\label{T2-TT:global}
\end{equation}

\PGRAPH{$\bullet$ TT format of $\mathbf{T}^{g,int}_3$}, $\ten{T}^{TT}_3$:
Following the formulation in Eqn.~\eqref{Convec:coef:Omega}, the convection operator in TT-matrix format is constructed as:
\begin{equation}
\begin{split}
    \ten{T}_3^{TT} 
=\Bigg (\sum_{m_1,m_2,m_3,m_4=0}^{1}  \big[&\bold{M}^{\text{g}}_{t_{m_4}}(\ten{I}_t,\ten{I}_t)\tenprod \bold{M}^{\text{g}}_{z_{m_3}}(\ten{I}_s,\ten{I}_s)\tenprod \bold{M}^{\text{g}}_{y_{m_2}}(\ten{I}_s,\ten{I}_s) \tenprod \bold{D}^{\text{g}}_{x_{m_1}}(\ten{I}_s,\ten{I}_s) 
\\
+&\bold{M}^{\text{g}}_{t_{m_4}}(\ten{I}_t,\ten{I}_t)\tenprod \bold{M}^{\text{g}}_{z_{m_3}}(\ten{I}_s,\ten{I}_s) \tenprod \bold{D}^{\text{g}}_{y_{m_2}}(\ten{I}_s,\ten{I}_s) \tenprod \bold{M}^{\text{g}}_{x_{m_1}}(\ten{I}_s,\ten{I}_s)  
\\
+&\bold{M}^{\text{g}}_{t_{m_4}}(\ten{I}_t,\ten{I}_t) \tenprod \bold{D}^{\text{g}}_{z_{m_3}}(\ten{I}_s,\ten{I}_s)\tenprod \bold{M}^{\text{g}}_{y_{m_2}}(\ten{I}_s,\ten{I}_s) \tenprod \bold{M}^{\text{g}}_{x_{m_1}}(\ten{I}_s,\ten{I}_s)\big]  \Bigg ) \ten{U}^{TT}.
\end{split}
\label{T3-TT:global}
\end{equation}

\PGRAPH{$\bullet$ TT format of $\mathbf{T}^{g,int}_4$}, $\ten{T}^{TT}_4$: Following the formulation in Eqn.~\eqref{T4:global}, the reaction term is constructed as:
\begin{equation}
\label{T4-TT:global}
    \ten{T}^{TT}_4:= \sum_{m_1,m_2,m_3,m_4=0}^{1}   \bigg(\bold{M}^{\text{g}}_{t_{m_4}}(\ten{I}_t,\ten{I}_t) \tenprod \bold{M}^{\text{g}}_{z_{m_3}}(\ten{I}_s,\ten{I}_s) \tenprod \bold{M}^{\text{g}}_{y_{m_2}}(\ten{I}_s,\ten{I}_s) \tenprod \bold{M}^{\text{g}}_{x_{m_1}}(\ten{I}_s,\ten{I}_s)\bigg) \ten{U}^{TT}.
\end{equation}
\PGRAPH{$\bullet$ TT format of $\mathbf{T}^{g,int}_5$}, $\ten{T}_5^{TT}$: From Eqn.~\eqref{T5:global}, the loading term is constructed as:
\begin{equation}
\label{T5-TT:global}
    \ten{T}^{TT}_5:= \bigg(\bold{M}_{I_T}(\ten{I}_t,\ten{I}_t) \tenprod \bold{M}_{\Omega_Z}(\ten{I}_s,\ten{I}_s) \tenprod \bold{M}_{\Omega_Y}(\ten{I}_s,\ten{I}_s) \tenprod \bold{M}_{\Omega_X}(\ten{I}_s,\ten{I}_s) \bigg) \ten{F}^{TT}.
\end{equation}
The TT format of loading term $\ten{F}^{TT}$ is approximated by the cross interpolation algorithm~\cite{Oseledets:2023}.

\PGRAPH{$\bullet$ TT format of the boundary term}, $\ten{F}^{Bd,TT}$: 
The boundary and initial conditions are enforced by the boundary term, which will map from all nodes to only interior nodes.
The detailed construction of $\ten{F}^{Bd,TT}$ is included in~\ref{appA}.

\subsection{Construction of QTT Format for Linear System Components}
\label{sec:QTT tensorization}
Given the banded structure of the component matrices, they can be further compressed through transformation into QTT format~\cite{kazeev2013multilevel}, we proceed with the construction of the QTT format.
The idea is to first convert the component matrices into QTT format and then linked together as in the case of TT format. This algorithm is described in~\cite[Algorithm 1]{truong2023tensor}. 

At this point, we have completed constructing the TT-format and QTT-format of the linear system in Eqn.~\eqref{eqn:TT/QTT linear system}.
To solve the TT/QTT linear system using tensor network optimization techniques, we have employed the MATLAB TT-Toolbox \cite{Oseledets:2023}.

\section{Numerical Experiments:}
\label{sec:numerical}
In this section we discuss a series of numerical experiments to assess the performance of our proposed tensor network space-time spectral element method, in both TT-format and QTT-format. All simulations are conducted using MATLAB software running on a Linux operating system equipped with a 2.1 GHz Intel Gold 6152 processor.
We compare the performance of TT and QTT solvers against the full grid solver, which performs calculations on full tensors. Further, we computed the error in $L^2([0,T] \times \Omega)$ norm for all numerical experiments. We have used linear polynomials for all experiments, and consequently, we have observed a quadratic order of convergence \cite[Figure~3]{gomez2024space}.  

\subsection{TT-Ranks of the Diffusion Operator}
In this experiment, we explore how the TT-ranks of the diffusion operator $\nabla \cdot \Big( \boldsymbol{\kappa}(x,y,z) \nabla u \Big)$ in TT format change with four different coefficient functions $\kappa(x,y,z)$. We choose to perform this investigation because it might be expected that the TT-ranks of the diffusion operator would not be particularly low, given the complex assembly process described in Equation~\eqref{T2-TT:global}. The chosen coefficient functions have varying TT-ranks to illustrate their impact on the TT-ranks of the diffusion operator. Moreover, these functions possess exact TT-ranks that are independent of grid sizes. The function $1/(1 + x + y + z)$, however, demonstrates more complex behavior, with TT-ranks that vary depending on the selected TT truncation tolerance, while still remaining independent of the grid size.

The result from Table~\ref{tab:TT-rank of diffusion operator}, shows that when the coefficient function is simply $\boldsymbol{\kappa}(x,y,z) = 1$, the TT-ranks of the global diffusion operator are surprisingly low at $[2, 2]$ across three discretization grid sizes. This demonstrates that despite the complex assembly process in SEM, the global diffusion operator itself possesses a low-rank structure, which however depends on the rank of the diffusion function and its rank.
In general, we observe that the TT-ranks of the diffusion operator increase at a rate that is at most twice the TT-ranks of the coefficient function $\boldsymbol{\kappa}(x,y,z)$. Moreover, the number of elements per dimension does not affect the TT-ranks, which will lead to better compression on a larger grid.
This result suggests a relationship between the complexity of the coefficient function and the rank of the resulting  global diffusion operator in the TT-format.

\begin{table}[httb]
    \centering
    \begin{tabular}{|c|c|c|c|c|c|}\hline
         \multirow{2}{*}{$\boldsymbol{\kappa}(x,y,z)$} & \multirow{2}{*}{TT-tol}  & \multirow{2}{*}{$\kappa^{TT}$ ranks}  & \multicolumn{3}{|c|}{TT-ranks of diffusion operator} \\\cline{4-6}
         & & & $N_q=17$ & $N_q=33$ & $N_q=65$\\\hline
        $1$ & 1e-12  & [1,1]  & [2,2]  & [2,2] & [2,2] \\\hline
        $1 + xyz$ & 1e-12  & [2,2]  & [4,4]  & [4,4] & [4,4] \\\hline   
        $1 + \cos(\pi(x+y))\cos(\pi z)$ & 1e-12  & [3,2]  & [6,4]  & [6,4] & [6,4] \\\hline
        $1/(1 + x+y+z)$ & 1e-6  & [5,5]  & [8,8]  & [8,8] & [8,8] \\\hline   
        $1/(1 + x+y+z)$ & 1e-12  & [9,9]  & [15,15]  & [15,15] & [15,15] \\\hline   
    \end{tabular}
    \caption{TT-ranks of diffusion operator. The result shows the connection between the TT-ranks of diffusion coefficient function $\boldsymbol{\kappa}(x,y,z)$, calculated via TT-cross interpolation with a truncation error TT-tol of the diffusion operator (the operator in the left-hand-side of Eq. \eqref{Poiason_Eq} in discrete form, and here the number of elements, $N_Q = N^{3}_q$, each element, $q_l$, has $N=2$ nodes.}
    \label{tab:TT-rank of diffusion operator}
\end{table}

\subsection{Three-Dimensional Poisson Equation}
Next, we demonstrate the performance of the TT and QTT solver on  the 3D Poisson equation,
\begin{equation}
  \begin{split}
 -\nabla \cdot \Big( \boldsymbol{\kappa}(x,y,z) \nabla u \Big)&= f(x,y,z) \ \text{in} \ \Omega,\\
    u&=0 \quad \text{on} \ \partial \Omega,\\
  \end{split}
  \label{Poiason_Eq}
\end{equation}
and with the manufactured solution $u(x,y,z)= \sin(\pi x)\sin( \pi y)\sin(\pi z)$ on the computational domain $[0,1] \times [0,1] \times [0,1]$, and $\boldsymbol{\kappa}(x,y,z) = 1 + \cos(\pi(x+y))\cos(\pi z)$. The loading term $f(x,y,z)$ and the boundary conditions are computed in the way to enforce the manufactured solution. We computed the solutions with three different solvers: full grid solver, TT solver and QTT solver.
\begin{figure}[httb]
    \centering
    \includegraphics[width=1.0\textwidth]{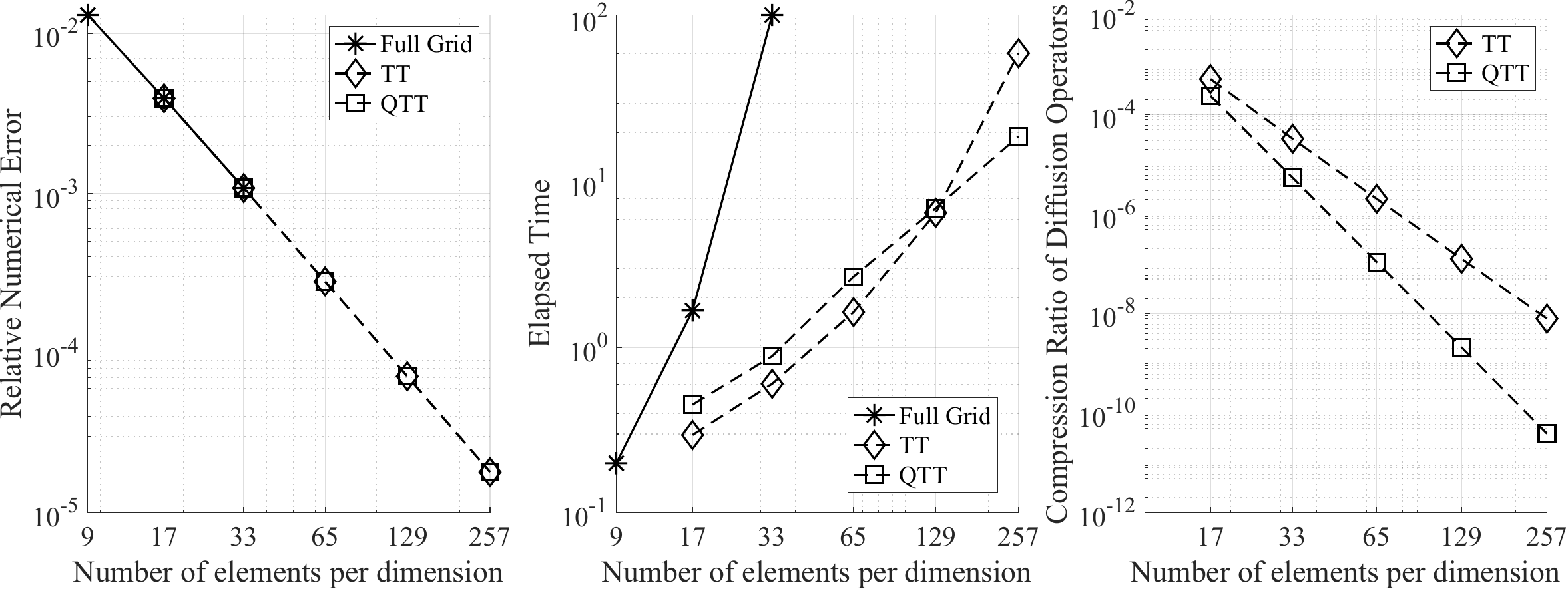}
    \caption{Performance of full-grid, TT and QTT solvers on 3D Poisson equation.}
    \label{fig:test1}
\end{figure}

The performance of three solvers are shown in Figure~\ref{fig:test1}, with convergence analysis in the left panel, the computational cost in the middle panel, and the compression ratio of TT and QTT format in the right panel.
Figure~\ref{fig:test1}-\textit{left} illustrates that both the TT and QTT formats maintain the same level of accuracy as the full grid solver up to a grid size of 33 elements per dimension. Beyond this point, the full grid solver requires more memory storage than is available. Both TT and QTT solvers were able to perform the simulation on finer grids, and successfully maintain the rate of convergence.

Figure~\ref{fig:test1}-\textit{middle} illustrates the computational time required by all three solvers. As expected, the full grid solver is the most computationally expensive, with a much steeper scaling compared to the TT and QTT solvers. When comparing the TT and QTT solvers, the TT solver is slightly faster for smaller grids but begins to show steeper scaling than the QTT solver at grid sizes of 129 and 257 elements per dimension.

Figure~\ref{fig:test1}-\textit{right} displays the compression ratio of the diffusion operators in TT and QTT formats across different grid sizes. The compression ratio is defined as the ratio of the number of elements in the full tensor to the number of elements in the TT or QTT format. The plot indicates that the QTT format becomes increasingly more efficient than the TT format as the grid size grows. For smaller grids, the compression achieved by the QTT format is similar to that of the TT format, which explains why the QTT solver is slightly slower than the TT solver at these grid sizes. This is because, at similar levels of compression, the QTT format involves calculations with more dimensions, making it more computationally demanding than the TT algorithm. However, as the compression advantage of the QTT format becomes more pronounced (about 2 order of magnitude lower at the grid size of 129), its computational efficiency starts to outweigh that of the TT format, reflecting in the reduced computational time.
%
\subsection{Three-Dimensional CDR Equation}
In this experiment, we investigate the performance of TT and QTT solvers on solving a 3D convection diffusion reaction equation with inhomogeneous boundary conditions.
\begin{equation}
  \begin{split}
  \frac{\partial u}{\partial t}  -\nabla \cdot \Big( \boldsymbol{\kappa}(t,\bold{x}) \nabla u \Big)+\vec{b} (t,\bold{x}) \cdot \nabla u 
   +c(t,\bx)u &= f(t,\bx) \ \text{in} \ \Omega \times [0,T],\\
    u&=g(t, \bold{x}) \quad \text{on} \ \partial \Omega \times [0,T],\\
   u(t=0, \bold{x})&=u_0(\bold{x}) \quad \text{in} \ \Omega,
  \end{split}
    \label{model:prob}
\end{equation}
where $\boldsymbol{\kappa}(t,\bold{x}) = 1 + cos(\pi x)cos(\pi y)cos(\pi z)$, $\vec{b} (\bold{x}) = [x,y,z]$, and $c(t,\bx) = e^{-(x+y+z)}$\\
The manufactured solution is $u(t,x,y,z) = sin(\pi(t+x+y+z))$ and the space time computational domain is $[0,1]^4$.

\begin{figure}[httb]
    \centering
    \includegraphics[width=1.0\textwidth]{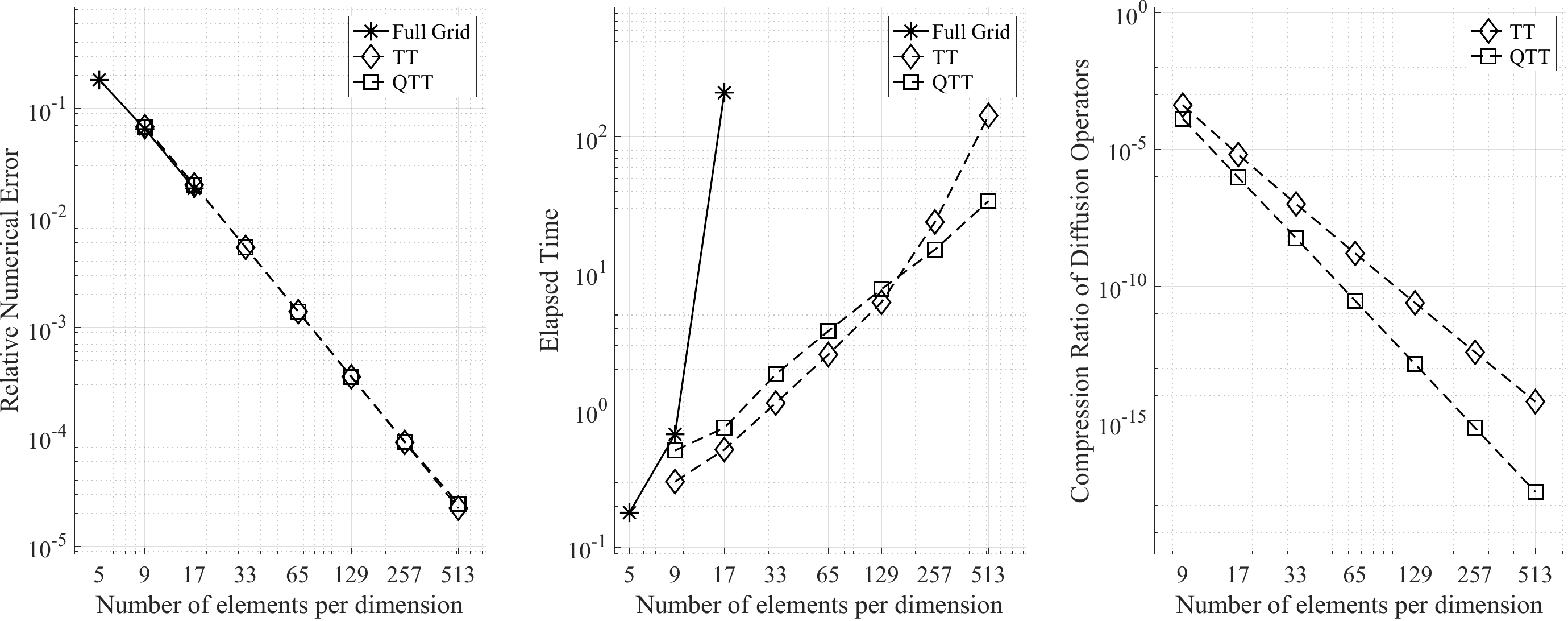}
    \caption{Performance of full-grid, TT and QTT solvers on 3D CDR equation with manufactured solution.}
    \label{fig:test2}
\end{figure}
The performance of three solvers (full grid, TT and QTT) is presented in Figure~\ref{fig:test2}, with convergence shown in the left panel, computational cost in the middle panel, and the compression ratio in the TT and QTT formats displayed in the right panel.

Figure~\ref{fig:test2}-\textit{left} demonstrates that both the TT and QTT formats maintain a similar level of accuracy as the full grid solver, consistent with the findings in the first experiment. Up to a grid size of 17 elements per dimension, all solvers perform comparably. However, beyond this point, the full grid solver becomes infeasible due to memory limitations. In contrast, both TT and QTT solvers continue to perform on finer grids, preserving the convergence rate even at higher resolutions.

Figure~\ref{fig:test2}-\textit{middle} shows the computational time required by the three solvers. The full grid solver is the most expensive one, exhibiting significantly steeper scaling with increasing grid size. When comparing the TT and QTT solvers, both TT and QTT solvers initially shows to be comparable for smaller grids. However, as grid sizes increase to 257 and 513 elements per dimension, the QTT solver's scaling becomes less steep, eventually overtaking the TT solver in computational efficiency.

Figure~\ref{fig:test2}-\textit{right} illustrates the behavior of the diffusion operators in TT and QTT formats across different grid sizes. The results show that the QTT format becomes increasingly more effective as the grid size increases. For example, at the grid size of 513 elements per dimension, the compression of QTT format is 4 order of magnitude lower compare to the compression of TT format.

Given that this is a 4D problem, the TT and QTT solvers has shown superior efficiency compared to the full grid solver, which potentially allows for much faster and more accurate simulations at much higher resolutions.

\subsection{Three-Dimensional CDR Equation with a Nonlinear Loading Term}

In this section, we present numerical results for a simplified semiconductor problem in where we replace drift-diffusion physics with the Poisson-Boltzmann equation approximation \cite{holst1993multilevel}. In this example, we benchmark the performance of TT and QTT solvers on a three-dimensional nonlinear equation
\begin{equation}
  \begin{split}
 \frac{\partial u}{\partial t} -\Delta u &=  u - u^3 + f(\bx) \ \text{in} \ \Omega \times [0,T],\\
    u&=0 \quad \text{on} \ \partial \Omega,\\
    u&(t=0, \bold{x})= u_0.
  \end{split}
\end{equation}
with the computational domain $[0,1]^4$, and the manufactured solution, \[u(t,x,y,z) = sin(\pi x)sin(\pi y)sin(\pi z)sin(\pi t) + sin(2\pi x)sin(2\pi y)sin(2\pi z)sin(2\pi t).\]
The output of discretization process is the nonlinear equation in TT format:
\[ \ten{A}^{TT}U^{TT} = \ten{M}^{TT}(\ten{U}^{TT}-(\ten{U}^{TT})^3 + \ten{F}^{TT}, \]
where $\ten{A}^{TT}U^{TT} = \ten{T}_1^{TT} + \ten{T}_2^{TT}$ and $\ten{M}^{TT} = \bold{M}_{I_T}(\ten{I}_t,\ten{I}_t) \tenprod \bold{M}_{\Omega_Z}(\ten{I}_s,\ten{I}_s) \tenprod \bold{M}_{\Omega_Y} \otimes \bold{M}_{\Omega_X}(\ten{I}_s,\ten{I}_s)$.
The QTT format of this equation can be achieved by using the process described above in Sec.~\ref{sec:QTT tensorization}.
\begin{figure}[httb]
    \centering
    \includegraphics[width=1.0\textwidth]{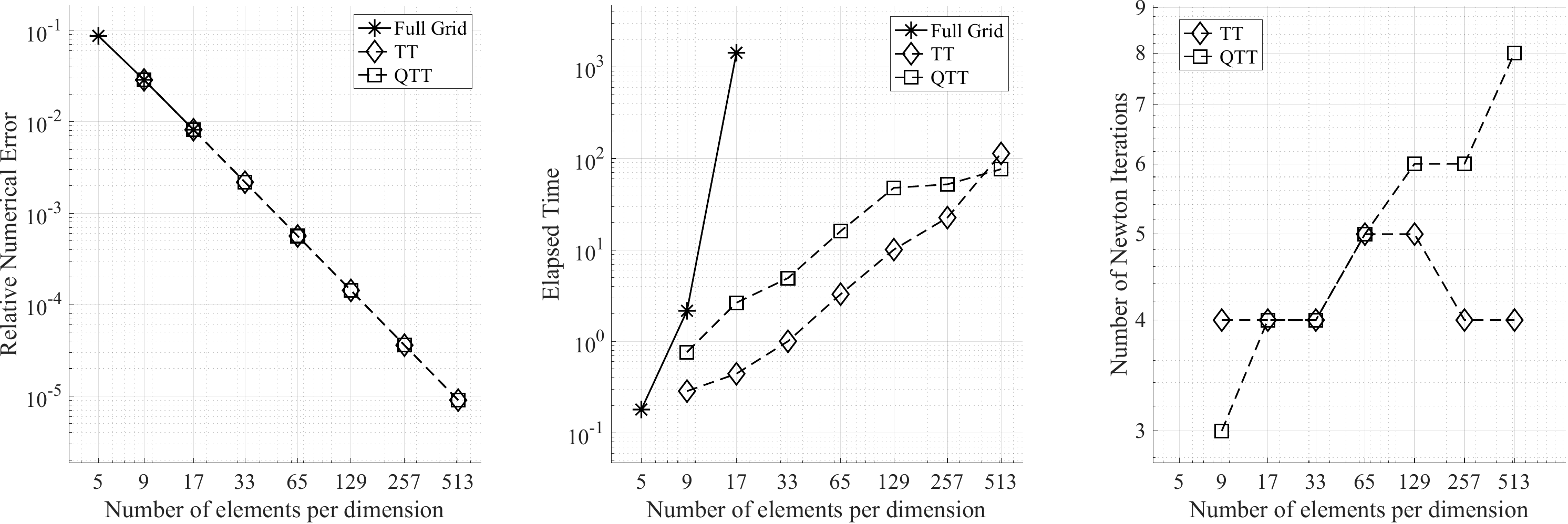}
    \caption{Performance of full-grid, TT and QTT solvers on 3D nonlinear equation with manufactured solution.}
    \label{fig:test3}
\end{figure}

The nonlinear TT/QTT equation is solved using the step truncation TT-Newton method developed in \cite[Algorithm1]{adak2024nonlinear}, with the loss function defined as:
\[\ten{L}(U^{TT}) = \ten{A}^{TT}U^{TT} - \ten{M}^{TT}(\ten{U}^{TT}-(\ten{U}^{TT})^3+\ten{F}^{TT} ).\]
The Jacobian of this matrix function is explicitly computable.
Figure~\ref{fig:test3} compares three solvers for this nonlinear problem, showing convergence (left panel), computational cost (middle panel), and the number of Newton iterations (right panel). The diffusion operator compression for TT and QTT solvers are not shown, as they are similar to the third experiment.

Figure~\ref{fig:test3}-\textit{left} demonstrates that both TT and QTT maintain accuracy comparable to the full grid solver up to 17 elements per dimension, beyond which the full grid solver becomes memory-infeasible. As expected, TT and QTT solvers efficiently handle finer grids and preserve convergence rates.

Figure~\ref{fig:test3}-\textit{middle} presents the computational time, where the full grid solver is the most expensive. The TT solver is initially faster for smaller grids, but QTT shows better scaling at larger grid sizes (257 and 513 elements per dimension).

However, as shown in Figure~\ref{fig:test3}-\textit{right}, at larger grid sizes (257 and 513), QTT solver requires more number of Newton iterations to converge.
For example, at the grid size of 513, QTT requires 8 iterations compared to 4 iterations of TT solver. This results in QTT being twice as fast in computational time per Newton iteration at this grid size. 
The results suggest that QTT and TT solvers may take distinct paths to convergence, which is intriguing and likely linked to the choice of truncation errors in these formats. This relationship is particularly important, as prior research has demonstrated that the selection of TT truncation errors significantly influences both the efficiency and accuracy of TT-format PDE solvers~\cite{rodgers2022adaptive,rodgers2023implicit,adak2024nonlinear,danis2024tensortrain,danis2024SW}. This observed difference in solver behavior certainly deserves more detailed investigation and thorough analysis.

\section{Conclusion}
\label{sec:conclusion}

In conclusion, we have developed a Petrov-Galerkin TT/QTT space-time formulation specifically tailored for the spectral element discretization of convection-diffusion-reaction equations. This approach includes a detailed construction of spectral element discrete operators in both TT and QTT formats, accommodating general coefficient functions. By leveraging the QTT format, we have effectively compressed the discrete operators while addressing the banded structure, leading to enhanced computational efficiency. Our numerical experiments demonstrate that our mixed approach not only achieves significant speedup but also drastically reduces memory usage, enabling high-resolution simulations to be conducted in a much shorter time. These advancements highlight the potential of our approach to tackle more complex and large-scale problems efficiently. Moving forward, future work will focus on extending this framework to address nonlinear partial differential equations and adapting it for problems in more complex domains, thereby broadening the scope and impact of this method in computational science and engineering.

We also plan to develop a discontinuous Petrov-Galerkin (DPG) formulation -- of the type introduced by Demkowicz and Gopalakrishnan \cite{DPG1,DPG2} -- by employing the so-called \emph{mixed formulation} of DPG (see Roberts et al. \cite[Appendix A.5.1]{RobertsVlasov1D1Vtimemarching}).  The mixed formulation will allow the application of TT formats by preserving the tensor-product structure of the basis functions; the mixed formulation avoids explicit computation of DPG's optimal test functions and computes an explicit error-representation function, which can be used to evaluate the accuracy of the solution.  In the TT context, the latter could enable an iterative refinement approach, whereby error estimates are used to determine placement of refined elements in the appropriate 1D mesh components.  DPG has previously been applied to \emph{steady} convection-diffusion-reaction problems \cite{bui2013unified}, and to space-time convection \cite{Broersen_Dahmen_Stevenson_18} and convection-reaction problems \cite{DemkowiczRobertsMatute}.  A DPG-based TT formulation of space-time convection-diffusion-reaction can be expected to provide excellent stability properties, even on coarse meshes and for challenging problems.

\section*{Data availability statement}

The data that support the findings of this research are available from the corresponding author upon reasonable request.

\section*{Declaration of competing interest}

The authors declare that they have no known competing financial interests or personal relationships that could have appeared to influence the work reported in this paper.

\section*{CRediT authorship contribution statement}
D. Adak, D. P. Truong, R. Vuchkov, S. De, D. DeSantis, N. V. Roberts, Kim O. Rasmussen, B. S. Alexandrov : Conceptualization,
Methodology, Writing-original draft, Review $\&$ Editing.


\section*{Acknowledgments}
The authors gratefully acknowledge the support of the Laboratory Directed Research and Development (LDRD) program of Los Alamos National Laboratory under projects number 20230067DR, 20240705ERil, and 20240782ER and ISTI Rapid Response 20248215CT-IST and Institutional Computing. Los Alamos National Laboratory is operated by Triad National Security, LLC, for the National Nuclear Security Administration of U.S. Department of Energy (Contract No. 89233218CNA000001).
\\
This research was sponsored by the U.S.~Department of Energy
Office of Science and the U.S.~Air Force Office of Scientific Research.
This article has been authored by employees of National Technology \&
Engineering Solutions of Sandia, LLC under Contract No.~DE-NA0003525 with
the U.S.\ Department of Energy (DOE). The employee owns all right, title
and interest in and to the article and is solely responsible for its
contents. The United States Government retains and the publisher, by
accepting the article for publication, acknowledges that the United States
Government retains a non-exclusive, paid-up, irrevocable, world-wide
license to publish or reproduce the published form of this article or allow
others to do so, for United States Government purposes. The DOE will
provide public access to these results of federally sponsored research in
accordance with the DOE Public Access Plan. Supported by the LDRD Program at Sandia National Laboratories. Sandia is managed and operated by
NTESS under DOE NNSA contract DE-NA0003525

\bibliographystyle{plain}
\bibliography{main}


\appendix
\section{Notations and Definitions}
\label{appdef}

\subsection{Kronecker product}
The Kronecker product $\bigotimes$ of matrix
$\mat{A}=(a_{ij})\in\mathbb{R}^{m_A \times n_A}$ and matrix
$\mat{B}=(b_{ij})\in\mathbb{R}^{m_B \times n_B}$ is the matrix
$\mat{A}\otimes \mat{B}$ of size
$N_{\mat{A}\otimes\mat{B}}=(m_Am_B)\times(n_An_B)$ defined as:
\begin{equation}
  \mat{A}\otimes\mat{B} = 
  \begin{bmatrix}
    a_{11}\mat{B} &a_{12}\mat{B} &\cdots & a_{1n_A}\mat{B}\\
    a_{21}\mat{B} &a_{22}\mat{B} &\cdots & a_{2n_A}\mat{B}\\
    \vdots & \vdots & \ddots & \vdots \\
    a_{m_A1}\mat{B} &a_{m_A2}\mat{B} &\cdots & a_{In_A}\mat{B}\\
  \end{bmatrix}.
  \label{def:kronecker}
\end{equation}
Equivalently, it holds that
$\big(\mat{A}\otimes\mat{B}\big)_{ij}=a_{i_Aj_A}b_{i_Bj_B}$, where
$i=i_B+(i_A-1)m_B$, $j=j_B+(j_A-1)m_B$, with $i_A=1,\ldots,m_A$,
$j_A=1,\ldots,n_A$, $i_B=1,\ldots,m_B$, and $j_B=1,\ldots,n_B$.

\subsection{The Tensor Product}
The tensor product of matrix $\mat{A}=(a_{ij})\in\mathbb{R}^{m_A \times n_A}$ and matrix
$\mat{B}=(b_{kl})\in\mathbb{R}^{m_B \times n_B}$ produces the
four-dimensional tensor of size $N_{\mat{A}\tenprod\mat{B}}=m_A \times
n_A \times m_B \times n_B$, with elements:
\begin{equation}
  \big(\mat{A}\tenprod\mat{B}\big)_{ijkl} = a_{ij}b_{kl},
  \label{def:outer_product:matrices}
\end{equation}
for
$i=1,2,\ldots,m_A$,
$j=1,2,\ldots,n_A$,
$k=1,2,\ldots,m_B$,
$l=1,2,\ldots,n_B$.

\section{TT Format Construction of the boundary term $\ten{F}^{bd,TT}$}
\label{appA}
The boundary term $\ten{F}^{bd,TT}$ is computed as \[\ten{F}^{bd,TT} = \ten{A}^{map,TT}\ten{G}^{bd,TT}\]
where $\ten{A}^{map,TT} = \ten{A}_t^{TT} + \ten{A}_d^{TT} + \ten{A}_c^{TT} + \ten{A}_r^{TT}.$

Next, we will show how to construct each component $\ten{A}_t^{TT},\ \ten{A}_d^{TT},\ \ten{A}_c^{TT},\ \ten{A}_r^{TT}.$ 
We redefine the sets of indices $\ten{I}_t = 2:N+1$ and $\ten{I}_s = 2:N$ where $N$ is the number of elements per dimension. Then the TT format of each component is constructed as:

\begin{equation}
  \ten{A}_t^{TT}:= \bold{D}_{{I_T}}(\ten{I}_t,:) \tenprod  \bold{M}_{\Omega_Z}(\ten{I}_s,:) \tenprod \bold{M}_{\Omega_Y}(\ten{I}_s,:) \tenprod \bold{M}_{\Omega_X}(\ten{I}_s,:).
\end{equation}

\begin{equation}
\begin{split}
    \ten{A}_d^{TT}:= \sum_{m_1,m_2,m_3,m_4=0}^{1}   \big[&\bold{M}^{\text{g}}_{t_{m_4}}(\ten{I}_t,:) \tenprod \bold{S}^{\text{g}}_{z_{m_3}}(\ten{I}_s,:) \tenprod \bold{M}^{\text{g}}_{y_{m_2}}(\ten{I}_s,:) \tenprod \bold{M}^{\text{g}}_{x_{m_1}}(\ten{I}_s,:)
 \\
 +&\bold{M}^{\text{g}}_{t_{m_4}}(\ten{I}_t,:) \tenprod \bold{M}^{\text{g}}_{z_{m_3}}(\ten{I}_s,:) \tenprod \bold{S}^{\text{g}}_{y_{m_2}}(\ten{I}_s,:) \tenprod \bold{M}^{\text{g}}_{x_{m_1}}(\ten{I}_s,:)\\
 +&\bold{M}^{\text{g}}_{t_{m_4}}(\ten{I}_t,:) \tenprod \bold{M}^{\text{g}}_{z_{m_3}}(\ten{I}_s,:) \tenprod \bold{M}^{\text{g}}_{y_{m_2}}(\ten{I}_s,:) \tenprod \bold{S}^{\text{g}}_{x_{m_1}}(\ten{I}_s,:)\big].
 \end{split}
\end{equation}

\begin{equation}
\begin{split}
    \ten{A}_c^{TT} 
=\sum_{m_1,m_2,m_3,m_4=0}^{1}  \big[&\bold{M}^{\text{g}}_{t_{m_4}}(\ten{I}_t,:)\tenprod \bold{M}^{\text{g}}_{z_{m_3}}(\ten{I}_s,:)\tenprod \bold{M}^{\text{g}}_{y_{m_2}}(\ten{I}_s,:) \tenprod \bold{D}^{\text{g}}_{x_{m_1}}(\ten{I}_s,:) 
\\
+&\bold{M}^{\text{g}}_{t_{m_4}}(\ten{I}_t,:)\tenprod \bold{M}^{\text{g}}_{z_{m_3}}(\ten{I}_s,:) \tenprod \bold{D}^{\text{g}}_{y_{m_2}}(\ten{I}_s,:) \tenprod \bold{M}^{\text{g}}_{x_{m_1}}(\ten{I}_s,:)  
\\
+&\bold{M}^{\text{g}}_{t_{m_4}}(\ten{I}_t,:) \tenprod \bold{D}^{\text{g}}_{z_{m_3}}(\ten{I}_s,:)\tenprod \bold{M}^{\text{g}}_{y_{m_2}}(\ten{I}_s,:) \tenprod \bold{M}^{\text{g}}_{x_{m_1}}(\ten{I}_s,:)\big].
\end{split}
\end{equation}

\begin{equation}
    \ten{A}_r^{TT}:= \sum_{m_1,m_2,m_3,m_4=0}^{1} \bold{M}^{\text{g}}_{t_{m_4}}(\ten{I}_t,:) \tenprod \bold{M}^{\text{g}}_{z_{m_3}}(\ten{I}_s,:) \tenprod \bold{M}^{\text{g}}_{y_{m_2}}(\ten{I}_s,:) \tenprod \bold{M}^{\text{g}}_{x_{m_1}}(\ten{I}_s,:).
\end{equation}

\begin{equation}
\ten{A}^{map,TT} = \ten{A}_t^{TT} + \ten{A}_d^{TT} + \ten{A}_c^{TT} + \ten{A}_r^{TT}.
\end{equation}

Next, we describe the construction of the tensor $\ten{G}^{bd,TT}$. The full grid tensor $G^{bd}$ is a $N_q \times (N-1) \times (N-1) \times (N-1)$ tensor, where only the initial and boundary elements are non-zero, with all other elements being zeros. The tensor $\ten{G}^{bd,TT}$ in TT format is then constructed using the cross interpolation technique.

Finally, the TT format of the boundary term is computed as $\ten{F}^{bd,TT} = \ten{A}^{map,TT}\ten{G}^{bd,TT}$.

\section{QTT Format Construction}
\label{appB}
\begin{equation}
\label{T1-QTT:global}
  \ten{T}_1^{QTT}:=   \Big ( \bold{D}^{QTT}_{{I_T}}(\ten{I}_t,\ten{I}_t) \tenprod  \bold{M}^{QTT}_{\Omega_Z}(\ten{I}_s,\ten{I}_s) \tenprod \bold{M}^{QTT}_{\Omega_Y}(\ten{I}_s,\ten{I}_s) \tenprod \bold{M}^{QTT}_{\Omega_X}(\ten{I}_s,\ten{I}_s) \Big) \ten{U}^{QTT},
\end{equation}

\begin{equation}
\begin{split}
    \ten{T}_2^{QTT}:= \Bigg ( &\sum_{m_1,m_2,m_3,m_4=0}^{1}   \big[\bold{M}^{\text{g},QTT}_{t_{m_4}}(\ten{I}_t,\ten{I}_t) \tenprod \bold{S}^{\text{g},QTT}_{z_{m_3}}(\ten{I}_s,\ten{I}_s) \tenprod \bold{M}^{\text{g},QTT}_{y_{m_2}}(\ten{I}_s,\ten{I}_s) \tenprod \bold{M}^{\text{g},QTT}_{x_{m_1}}(\ten{I}_s,\ten{I}_s)
 \\
 +&\bold{M}^{\text{g},QTT}_{t_{m_4}}(\ten{I}_t,\ten{I}_t) \tenprod \bold{M}^{\text{g},QTT}_{z_{m_3}}(\ten{I}_s,\ten{I}_s) \tenprod \bold{S}^{\text{g},QTT}_{y_{m_2}}(\ten{I}_s,\ten{I}_s) \tenprod \bold{M}^{\text{g},QTT}_{x_{m_1}}(\ten{I}_s,\ten{I}_s)\\
 +&\bold{M}^{\text{g},QTT}_{t_{m_4}}(\ten{I}_t,\ten{I}_t) \tenprod \bold{M}^{\text{g},QTT}_{z_{m_3}}(\ten{I}_s,\ten{I}_s) \tenprod \bold{M}^{\text{g},QTT}_{y_{m_2}}(\ten{I}_s,\ten{I}_s) \tenprod \bold{S}^{\text{g},QTT}_{x_{m_1}}(\ten{I}_s,\ten{I}_s)\big] \Bigg ) \ten{U}^{QTT}.
 \end{split}
\label{T2-QTT:global}
\end{equation}

\begin{equation}
\begin{split}
    \ten{T}_3^{QTT} 
=\Bigg (&\sum_{m_1,m_2,m_3,m_4=0}^{1}  \big[\bold{M}^{\text{g},QTT}_{t_{m_4}}(\ten{I}_t,\ten{I}_t)\tenprod \bold{M}^{\text{g},QTT}_{z_{m_3}}(\ten{I}_s,\ten{I}_s)\tenprod \bold{M}^{\text{g},QTT}_{y_{m_2}}(\ten{I}_s,\ten{I}_s) \tenprod \bold{D}^{\text{g},QTT}_{x_{m_1}}(\ten{I}_s,\ten{I}_s) 
\\
+&\bold{M}^{\text{g},QTT}_{t_{m_4}}(\ten{I}_t,\ten{I}_t)\tenprod \bold{M}^{\text{g},QTT}_{z_{m_3}}(\ten{I}_s,\ten{I}_s) \tenprod \bold{D}^{\text{g},QTT}_{y_{m_2}}(\ten{I}_s,\ten{I}_s) \tenprod \bold{M}^{\text{g},QTT}_{x_{m_1}}(\ten{I}_s,\ten{I}_s)  
\\
+&\bold{M}^{\text{g},QTT}_{t_{m_4}}(\ten{I}_t,\ten{I}_t) \tenprod \bold{D}^{\text{g},QTT}_{z_{m_3}}(\ten{I}_s,\ten{I}_s)\tenprod \bold{M}^{\text{g},QTT}_{y_{m_2}}(\ten{I}_s,\ten{I}_s) \tenprod \bold{M}^{\text{g},QTT}_{x_{m_1}}(\ten{I}_s,\ten{I}_s)\big]  \Bigg ) \ten{U}^{QTT},
\end{split}
\label{T3-QTT:global}
\end{equation}

\begin{equation}
\label{T4-QTT:global}
\begin{split}
\ten{T}^{QTT}_4:= \sum_{m_1,m_2,m_3,m_4=0}^{1}   \bigg(&\bold{M}^{\text{g},QTT}_{t_{m_4}}(\ten{I}_t,\ten{I}_t) \tenprod \bold{M}^{\text{g},QTT}_{z_{m_3}}(\ten{I}_s,\ten{I}_s)... \\ &\tenprod \bold{M}^{\text{g},QTT}_{y_{m_2}}(\ten{I}_s,\ten{I}_s) \tenprod \bold{M}^{\text{g},QTT}_{x_{m_1}}(\ten{I}_s,\ten{I}_s)\bigg) \ten{U}^{QTT}.
\end{split}
\end{equation}

\begin{equation}
\label{T5-QTT:global}
    \ten{T}^{QTT}_5:= \bigg(\bold{M}^{QTT}_{I_T}(\ten{I}_t,\ten{I}_t) \tenprod \bold{M}^{QTT}_{\Omega_Z}(\ten{I}_s,\ten{I}_s) \tenprod \bold{M}^{QTT}_{\Omega_Y} \otimes \bold{M}^{QTT}_{\Omega_X}(\ten{I}_s,\ten{I}_s) \bigg) \ten{F}^{QTT},
\end{equation}
The QTT format of loading term $\ten{F}^{QTT}$ is approximated by first approximate the $\ten{F}^{TT}$ using the cross interpolation algorithm, the convert the $\ten{F}^{TT}$ into the $\ten{F}^{QTT}$. The QTT format of the boundary term  $\ten{F}^{Bd,QTT}$ can be similarly constructed using~\ref{appA}.

\end{document}